\documentclass[11pt]{amsart}
\usepackage{amsfonts,newlfont,latexsym}
\usepackage{euscript}
\usepackage{amssymb}
\usepackage{amsmath}
\usepackage{color}
\usepackage{hyperref}
\definecolor{orange}{rgb}{1,0.5,0}

\numberwithin{equation}{section}

      \newtheorem{Theorem}{Theorem}[section]
       \newtheorem{theorem}[Theorem]{Theorem}
      \newtheorem{lemma}[Theorem]{Lemma}
      \newtheorem{Lemma}[Theorem]{Lemma}
      
      \newtheorem{proposition}[Theorem]{Proposition}
      
      \newtheorem{remark}[Theorem]{Remark}

\newcommand{\R}{\mathbb R}

\newcommand{\C}{\mathbb C}
\newcommand{\Z}{\mathbb Z}

\newcommand{\Q}{\mathbb Q}

\def\Proof{{\em Proof}\,: }
\def\proof{{\em Proof}\,: }
\def\QED{~\hfill~ $\diamond$ \vspace{7mm}}

\def \C {{\mathbb C }}

\def \R{{\mathbb R}}
\def \Z{{\mathbb Z}}
\def \N{{\mathbb N}}

\def \Q{{\mathbb Q}}

\def\holder{H\"{o}lder }

\definecolor{rjs}{rgb}{.1,.4,.7}

\begin{document}

\author[Alexander Gorodnik and Ralf Spatzier]
{Alexander Gorodnik$^\ast$ and  Ralf Spatzier$^{\ast \ast}$}

\title[Mixing Properties of Commuting Nilmanifold Automorphisms]
{Mixing Properties of Commuting Nilmanifold Automorphisms}

\thanks{$^\ast$ Supported in part by EPSRC grant EP/H000091/1 and ERC grant 239606}
\thanks{$^{\ast \ast}$ Supported in part by NSF grant DMS-0906085}

 \address{School of Mathematics, University of Bristol, Howard House, Queens Ave., Clifton, Bristol, BS8 1SD, U.K.}

\email{a.gorodnik@bristol.ac.uk}

 \address{Department of Mathematics, University of Michigan, Ann Arbor, MI 48109, USA.}

\email{spatzier@umich.edu}

\begin{abstract}
We study mixing properties of commutative groups of automorphisms acting on compact nilmanifolds.
Assuming that every nontrivial element acts ergodically, we prove that such actions are
mixing of all orders.  We further show exponential 2-mixing and 3-mixing. 
As an application we prove smooth cocycle rigidity for higher-rank abelian groups of nilmanifold automorphisms.
\end{abstract}

\maketitle

\section{Introduction}

Given a measure-preserving action of a (discrete) group $\Gamma$ on a probability space $(X,\mu)$,
we say that this action is \emph{$(s+1)$-mixing} if for every $f_0,\ldots, f_s\in L^\infty(X)$
and $\gamma_0,\ldots,\gamma_s\in\Gamma$,
\begin{equation}\label{eq:multi}
\int_X \left( \prod_{i=0}^s f_i(\gamma_ix)\right)\, d\mu(x)\longrightarrow \prod_{i=0}^s \left(\int_X f_i\, d\mu\right)
\end{equation}
as $\gamma_{i_1}\gamma_{i_2}^{-1}\to \infty$ for all $i_1\ne i_2$.
In particular, 2-mixing corresponds to the usual notion of mixing.
It was discovered by F.~Ledrappier \cite{Ledr} that 2-mixing does not imply 3-mixing in general.
In this paper we will be interested in mixing of higher order.
This property  is a very widespread phenomenon for one-parameter actions.
It is known to hold for many transformations satisfying some hyperbolicity assumptions.
On the other hand,
very little is known about higher order mixing for actions of large groups.
We are only aware of two general families of actions of large groups on manifolds
where the multiple mixing has been established ---
$\Z^l$-actions by automorphisms on compact abelian groups
and actions of simple Lie groups.
K. Schmidt and T. Ward \cite{sw} proved that 2-mixing $\Z^l$-action by automorphisms 
on compact connected abelian groups are mixing of all orders,
and S.~Mozes \cite{Mozes} established mixing of all orders for ergodic actions of simple Lie groups. 

In this paper we investigate mixing properties of $\Z^l$-actions by
automorphisms on compact nilmanifolds. We prove that for such actions, 2-mixing implies 
mixing of all orders and establish quantitative estimates for 2-mixing and 3-mixing.

\subsection{Main results}
Let $G$ be a simply connected nilpotent group and $\Lambda$
a discrete cocompact subgroup. We call the space $X=G/\Lambda$ a {\em compact nilmanifold}.
We denote by $\hbox{Aut}(X)$ the group of continuous automorphisms $\alpha$ of $G$
such that $\alpha(\Lambda)=\Lambda$. Then $\hbox{Aut}(X)$ naturally acts on $X$
and preserves the Haar probability measure $\mu$ on $X$.

Our first main result concerns exponential 3-mixing.
In order to obtain any quantitative estimate in (\ref{eq:multi}), it is necessary to work in a class
of sufficiently regular functions. We denote by $C^\theta(X)$ the space of
H\"older functions with exponent $\theta$, defined with respect to a Riemannian metric on $X$.

\begin{Theorem} \label{th:3_mixing} 
Let $\alpha:\Z^l\to \hbox{Aut}(X)$ be an action on a compact nilmanifold $X$
such that every $\alpha(z)$, $z\ne 0$, is ergodic.
Then there exists $\eta=\eta(\theta)>0$ such that for every $f_0,f_1,f_2\in C^\theta(X)$ and $z_0,z_1,z_2\in \mathbb{Z}^l$,
\begin{align*}
&\int_X f_0(\alpha(z_0)x)f_1(\alpha(z_1)x)f_2(\alpha(z_2)x)\,d\mu(x)\\
=& \left(\int_X f_0\, d\mu\right)\left(\int_X f_1\, d\mu\right)\left(\int_X f_2\, d\mu\right)
+O\left(N(z_0,z_1,z_2)^{-\eta}\|f_0\|_{C^\theta}\|f_1\|_{C^\theta}\|f_2\|_{C^\theta}\right),
\end{align*}
where $N(z_0,z_1,z_2)=\exp(\min_{i\ne j} \|z_i-z_j\|)$.
\end{Theorem} 

We note that this result is new even for the case of toral automorphisms.
Previously, quantitative 2-mixing  was established for toral automorphisms in \cite{Lind}
and for automorphisms of more general compact abelian groups in \cite{mw}.
Mixing of all orders for ergodic commuting toral automorphisms was established in \cite{sw}.
The argument in \cite{sw} relies on finiteness of the number of nondegenerate solutions of $S$-unit equations established
in \cite{Schlickewei}.
Although there are explicit estimates on the number of such solutions,
these estimates are not sufficient to derive any quantitative estimate for 3-mixing
because it is also essential to know
how the sets of solutions depend on the coefficients.
In order to prove Theorem \ref{th:3_mixing}, we use more delicate Diophantine estimates 
for linear forms in logarithms of algebraic numbers established in \cite{wald}
(cf. Proposition \ref{p:waldschmidt} below).

\vspace{0.2cm}

We also prove mixing of all orders:

\begin{Theorem} \label{th:multmixing}
Let $\alpha:\Z^l\to \hbox{Aut}(X)$ be an action on a compact nilmanifold $X$
such that every $\alpha(z)$, $z\ne 0$, is ergodic.
Then for every $f_0,\ldots,f_s\in L^\infty(X)$ and $z_0,\ldots,z_s\in \Z^l$, 
$$
\int_X \left( \prod_{i=0}^s f_i(\alpha (z_i)x)\right)\, d\mu(x)
= \prod_{i=0}^s \left(\int_X f_i\, d\mu\right)+o(1)
$$
as $\min_{i\ne j} \|z_i-z_j\|\to \infty$. Moreover, the convergence is uniform over families of  H\"older
functions $f_0,\ldots,f_s$ such that $\|f_0\|_{C^\theta},\ldots, \|f_s\|_{C^\theta}\ll 1$.
\end{Theorem}

This theorem extends the main result of \cite{sw} to general nilmanifolds.
The proof in \cite{sw} utilises  abelian Fourier analysis 
and properties of solutions of $S$-unit equations.
Our approach is based on the study of distribution of images of polynomial maps in $X$.
This reduces the proof to the investigation of certain Diophantine inequalities
which are analysed using W.~Schmidt's Subspace Theorem.
In order to prove an effective version of Theorem \ref{th:multmixing}, one would
need to estimate the size of nondegenerate solutions of these Diophantine
inequalities in terms of complexities of coefficients
(cf. Proposition \ref{p:subspace} below). However, this seems to be far out of reach
of available techniques when $s>2$.

\vspace{0.2cm}

Finally, we discuss the problem of exponential mixing for shapes in $\hbox{Aut}(X)$.
This notion was introduced by K.~Schmidt in \cite{Schmidt} in order to better understand
 Ledrappier's examples \cite{Ledr} which are not mixing of higher order.
A {\em shape} in $\hbox{Aut}(X)$ is a collection of  elements $\alpha_0, \ldots, \alpha_s \in \hbox{Aut}(X)$.
We says that the shape is {\it mixing} if for every $f_0,\ldots, f_s\in L^\infty(X)$,
$$
\int_X \left( \prod_{i=0}^s f_i(\alpha_i^n x)\right)\, d\mu(x)\longrightarrow \prod_{i=0}^s \left(\int_X
  f_i\, d\mu\right) 
$$
as $n\to \infty$. This property 
has been extensively studied in the context of commuting automorphisms of compact abelian groups
(see, for instance, \cite{Einsiedler-Ward}, \cite[Ch.~VIII]{Sch_book}, \cite{ward1}, \cite{ward2}).

We establish quantitative mixing for commuting Anosov shapes.
We say that the shape $\alpha_0,\ldots,\alpha_s$ 
is {\it Anosov} if $\alpha_{i_i}\alpha_{i_2}^{-1}$ is an Anosov map
for all $i_1\ne i_2$.

\begin{theorem}\label{th:shape}
Let $X$ be a compact nilmanifold and $\alpha_0,\ldots,\alpha_s\in \hbox{Aut}(X)$
a commuting Anosov shape. Then there exists $\rho=\rho(\theta)\in (0,1)$ such that
for every $f_0,\ldots, f_s\in C^\theta(X)$ and $n\in \N$,
\begin{equation}  \label{eq:shape}
 \int_X \left( \prod_{i=0}^s f_i(\alpha_i^n x)\right)\, d\mu(x) = \prod_{i=0}^s \left(\int_X f_i\, d\mu\right)+O\left( \rho^{  n   }            \prod_{i=0}^s \|f_i\|_{C^\theta} \right).
\end{equation}
\end{theorem}

\vspace{0.2cm}

\subsection{ Applications to rigidity}
The exponential mixing property played an important role in the program of 
classification of smooth Anosov higher-rank
$\Z^l$-actions on compact manifolds.
It is expected that all such actions can be built from actions by automorphisms on 
nilmanifolds. 
D. Fisher, B. Kalinin and R. Spatzier in \cite{FKS} applied the exponential 2-mixing property to extend their 
results for Anosov $\mathbb{Z}^l$-actions on tori to actions on nilmanifolds:

\begin{theorem}[D.~Fisher, B. Kalinin, R. Spatzier] \label{th:fks}
Let $\alpha$ be a $C^\infty$-action of $\Z^l$, $l\ge 2$, on a compact nilmanifold $X$
and let $\rho:\Z^l\to\hbox{Aut}(X)$ be the map induced by the action of $\alpha(\Z^l)$ on the 
fundamental group of $X$. Assume that there is a $\Z^2$ subgroup of $\Z^l$ such that $\rho(z)$
is ergodic for every nonzero $z\in \Z^2$, and there is an Anosov element for $\alpha$ in each Weyl
chamber of $\rho$. Then $\alpha$ is $C^\infty$-isomorphic to $\rho$.  
\end{theorem}

In fact, this application to global rigidity 
was our original motivation to establish the
exponential mixing property for nilmanifold automorphisms.  

Recently, F. Rodriguez-Hertz and Z. Wang \cite{rw} generalised Theorem \ref{th:fks} and
established a global rigidity result 
using only  existence of a single Anosov element.   
Again, they crucially use the exponential mixing property,
and reduce the problem to the prior result by showing
existence of many Anosov elements.

\vspace{0.2cm}

We also use the exponential mixing property to establish cocycle rigidity for 
higher-rank $\mathbb{Z}^l$-actions by automorphisms of nilmanifolds, extending the results
of A.~Katok and R.~Spatzier \cite{KS,KS2}. 
A {\it $C^{\infty}$-cocycle}  is a $C^{\infty}$-map $c:  \Z^l \times X \mapsto \R$ satisfying the identity 
\[  c (z_1+z_2, x) = c (z_1, z_2 \cdot x ) + c (z_2,x) \quad \hbox{for $z_1,z_2\in\Z^l$ and $x\in X$.}\]
Two cocycles $c_1$ and  $c_2$ are called {\em smoothly cohomologous} if there exists $b\in
C^\infty(X)$ such that
$$
c_1(z,x)=c_2(z,x)+b(z\cdot x)-b(x)\quad \hbox{for $z\in \Z^l$ and $x\in X$.}
$$
We call a  cocycle {\em constant} if it does not depend on $x \in X$.  
We prove that cocycles over genuine higher-rank actions by automorphisms on nilmanifolds
are smoothly cohomologous to constant cocycles.  This phenomenon 
 was 
first discovered by A. Katok and R. Spatzier in \cite{KS} for certain  higher-rank Anosov actions. Using
our methods, we generalise this cocycle rigidity theorem to actions by  automorphisms on nilmanifolds. 
We emphasise that the action in the following theorem  need not be Anosov. 

\begin{Theorem}\label{th:cocyle}
Let $\alpha:\Z^l\to \hbox{Aut}(X)$ be an action on a compact nilmanifold $X$.
Assume that there is a $\Z^2$ subgroup of $\Z^l$ such that $\alpha(z)$
is ergodic for every nonzero $z\in \Z^2$.
Then every smooth $\R$-valued cocycle is smoothly cohomologous to a constant cocyle.
\end{Theorem}

 For certain actions by partially hyperbolic left translations on homogeneous spaces $G/\Gamma$, 
where $G$ is a semisimple Lie groups and $\Gamma$ is a lattice in $G$, 
a similar theorem was
 proved by  D. Damjanovic and A. Katok
 \cite{DamKat05,DamKat10,Dam07} and Z. J. Wang 
\cite{Wang10}.  We note that these authors also prove \holder  versions of this result which are not amenable to our techniques. 
Furthermore, cocycle rigidity results are proven for small perturbations of these actions on $G/\Gamma$ in
\cite{DamKat10,Wang10}.  Again we cannot obtain these results by our methods.

\subsection*{ Acknowledgements}
A.G. would like to thank the University of Michigan for hospitality
during his visit when the work on this project had started.
R.S. thanks the University of Bristol for hospitality and support during this work.

\section{Exponential 3-mixing}\label{sec:3_m}

In this section we prove Theorem \ref{th:3_mixing}.
We start by setting up basic notation, which will be also used in subsequent sections.
Then in Sections \ref{sec:lyapunov}--\ref{Green-Tao dichotomy}, we collect some auxiliary estimates.
The proof of Theorem \ref{th:3_mixing} is divided into two parts.
We first give a proof under an irreducibility condition in Section \ref{sec:3_irred}, and then
in Section \ref{sec:3_gen} prove Theorem \ref{th:3_mixing} in general using an inductive argument.

We note that if the reader is only interested in exponential 2-mixing,
then the results of Section \ref{sec:dioph} are not needed, and in Section \ref{sec:3_irred},
one only needs to consider Case 1. This makes the proof much simpler.

\subsection{Notation}\label{sec:notation}
Let $G$ be a connected simply connected nilpotent Lie group, $\Lambda$ a discrete cocompact subgroup,
and $X=G/\Lambda$ the corresponding nilmanifold equipped with the invariant probability measure $\mu$. 
We fix a a right invariant Riemannian metric $d$ on $G$ which also defines a Riemannian metric on $X$.
Let $\mathcal{L}(G)$ be the Lie algebra of $G$ and $\exp:\mathcal{L}(G)\to G$ the exponential map. 
The lattice subgroup $\Lambda$ defines a rational structure on $\mathcal{L}(G)$.
For a field $K\supset \mathbb{Q}$,
we denote by $\mathcal{L}(G)_{K}$ the corresponding Lie algebra over $K$. 
Let $\pi:G\to G/G'$ denote the factor map, where $G'$ is the commutator subgroup.
We also have the corresponding map $D\pi:\mathcal{L}(G)\to \mathcal{L}(G/G')$.
We fix an identification $G/G'\simeq \mathcal{L}(G/G')\simeq \mathbb{R}^d$ that respects the rational structures.

Every automorphism $\beta$ of $G$ defines a Lie-algebra automorphism $D\beta:\mathcal{L}(G)\to
\mathcal{L}(G)$ such that $\beta\circ \exp=\exp\circ D\beta$.
If $\beta(\Lambda)= \Lambda$, then $D\beta$ preserves the
rational structure of $\mathcal{L}(G)$ defined by $\Lambda$.
In particular, given an action  $\alpha:\mathbb{Z}^l\to \hbox{Aut}(X)$ on the nilmanifold $X=G/\Lambda$,
we obtain a homomorphism $D\alpha:\mathbb{Z}^l\to \hbox{GL}(\mathcal{L}(G)_\Q)$.

For a multiplicative character $\chi:\mathbb{Z}^l\to \mathbb{C}^\times$, we set
$$
\mathcal{L}_\chi:=\{u\in \mathcal{L}(G)\otimes \mathbb{C}:\, D\alpha(z)\, u= \chi(z)u
\hbox{ for $z\in\mathbb{Z}^l$}\}.
$$
Let $\mathcal{X}(\alpha)$ denote the set of characters $\chi$ appearing in  the action $D\alpha$ on
$\mathcal{L}(G)$, and $\mathcal{X}'(\alpha)\subset \mathcal{X}(\alpha)$ be the set of
characters appearing in the action on $\mathcal{L}(G)/\mathcal{L}(G)'$.

\subsection{Estimates on Lyapunov exponents}\label{sec:lyapunov}
Since $\alpha(\Z^l)$ preserves the rational structure on $\mathcal{L}(G)$
defined by the lattice $\Lambda$,
it follows that each character
$\chi$ in $\mathcal{X}(\alpha)$ is of the form $\chi(z)=u_1^{z_1}\cdots u_l^{z_l}$
where $u_i$'s are algebraic numbers. The Galois group $\hbox{Gal}(\bar \Q/\Q)$
naturally acts on $\mathcal{X}(\alpha)$ and $\mathcal{X}'(\alpha)$. 
Let $\mathcal{X}_0\subset \mathcal{X}'(\alpha)$ be one of the Galois orbits.

\begin{Lemma}\label{l:expansion}
Suppose that every $\alpha(z)$, $z\ne 0$, acts ergodically on $X$.
Then there exists $c>0$ such that
$$
\max_{\chi\in \mathcal{X}_0} |\chi(z)|\ge \exp(c\|z\|)\quad \hbox{for all $z\in\Z^l$.}
$$
\end{Lemma}
 
\proof
By \cite[5.4.13]{CG}, $\Lambda G'/G'$ is a lattice in $G/G'$, and the action $\alpha$
defines the action on the torus $T:=G/(\Lambda G')$ by linear automorphisms.
Let $V$ be the subspace of $\mathcal{L}(G)/\mathcal{L}(G)'$ spanned by
the $\chi$-eigenspaces in with $\chi\in\mathcal{X}_0$.
Clearly, this subspace is invariant under $\alpha(\mathbb{Z}^l)$
	and is defined over $\mathbb{Q}$. Hence, it defines an $\alpha$-invariant subtorus 
	 $T_{\mathcal{X}_0}$     
of $T$. Since $\alpha(z)|_T$ is ergodic when $z\ne 0$,
it follows that the corresponding linear map has no roots of unity as eigenvalues.
This implies that  
$\alpha(z)|_{T_{\mathcal{X}_0}}$
 is also ergodic.


Consider a linear map $\ell: \mathbb{R}^l\to \mathbb{R}^{|\mathcal{X}_0|}$ which is defined
for $z\in \mathbb{Z}^l$ by
$$
\ell(z):=\left(\log |\chi(z)|:\, \chi\in \mathcal{X}_0\right)
$$
and extended to $\mathbb{R}^l$ by linearity.
Since for every $z\in \mathbb{Z}^l\backslash \{0\}$, the automorphism $\alpha(z)$ acting on $T_{\mathcal{X}_0}$
is ergodic, we have $\ell(z)\ne 0$ by \cite[Lemma 3.2]{GS}.
Hence, $\ell|_{\mathbb{Z}^l}$ is injective.

We also claim that $\ell(\mathbb{Z}^l)$ is discrete.
 We consider the embedding
$\mathbb{Z}^l \to \hbox{GL}(V)$ defined by $\alpha$.
Since $\alpha(\mathbb{Z}^l)$ preserves the integral lattice in $V$ corresponding to the torus $T_{\mathcal{X}_0}$,
it follows that the image of this embedding is discrete.
In other words, the subset
$$
\{(\chi(z):\, \chi\in \mathcal{X}_0):\, z\in\mathbb{Z}^l\}
$$
of $(\mathbb{C}^\times)^{|\mathcal{X}_0|}$ is discrete.
Since the kernel of the natural homomorphism $(\mathbb{C}^\times)^{|\mathcal{X}_0|}\to 
\mathbb{R}^{|\mathcal{X}_0|}$ defined by $s\mapsto \log |s|$ is compact,
this implies that $\ell(\mathbb{Z}^l)$ is discrete, as claimed.
Since $\ell(\mathbb{Z}^l)$ is discrete and has rank $l$,
it follows that  the space $\ell(\mathbb{R}^l)$ has dimension $l$,
and, in particular, the map $\ell$ is injective. Therefore, 
by compactness, there exists $c>0$
such that for every $z\in\mathbb{R}^l$, we have
$$
\max\{\log |\chi(z)|:\, \chi\in \mathcal{X}_0\}
\ge c\|z\|. 
$$
This implies the lemma.
\QED

Lemma \ref{l:expansion} shows that in Theorem \ref{th:3_mixing} we may replace $N(z_0,z_1,z_2)$
by 
$$
\exp\left(\min_{i\ne j} \max_{\chi\in \mathcal{X}_0} |\chi(z_i-z_j)|\right).
$$

\subsection{Diophantine estimates}\label{sec:dioph}

Recall that the (absolute) height of an algebraic number $u$ is defined by
$$
\hbox{H}(u)=\left(\prod_v \max(1, |u|_v)\right)^{1/[\Q(u):\Q]},
$$
where $|\cdot|_v$ denote suitably normalised absolute values of the field $\Q(u)$.
When $u$ is an algebraic integer, the height can be computed as
$$
\hbox{H}(u)=\left(\prod_i \max(1, |u_i|)\right)^{1/[\Q(u):\Q]},
$$
where $u_i$ denote all the Galois conjugates of $u$.

The following result is deduced from the work of M.~Waldschmidt \cite[Cor.~10.1]{wald}.

\begin{proposition}\label{p:waldschmidt}
Let $u_1,\ldots u_l, u\in\C$ be algebraic numbers and $z=(z_1,\ldots,z_l)\in\Z^l$.
Then there exists $c_1,c_2,c_3>1$, depending on $u_1,\ldots,u_l$ and $[\Q(u):\Q]$,
such that assuming that
\begin{equation}\label{eq:cond}
\|z\|\ge \log (c_2\textrm{\rm H}(u)),
\end{equation}
and 
$$
u_1^{z_1}\cdots u_l^{z_l}u\ne 1,
$$
we have the estimate
\begin{equation}\label{eq:ineq_main}
|u_1^{z_1}\cdots u_l^{z_l}u-1|\ge \exp\left(-c_1\log (c_2 \hbox{\rm H}(u))\log\left(\frac{c_3\|z\|}{\log
      (c_2\hbox{\rm H}(u))}\right)\right).
\end{equation}
\end{proposition}

Surprisingly, it turns out that the term $\log(c_2\hbox{\rm H}(u))$ in the denominator 
is essential to establish exponential 3-mixing (cf. (\ref{eq:e2})--(\ref{eq:kappa4}) below).

\vspace{0.2cm}

\proof
We note that since $\hbox{\rm H}(u)\ge 1$ and (\ref{eq:cond}) holds,
the right hand side of (\ref{eq:ineq_main}) is bounded from above by
$$
\exp(-c_1\log (c_2) \log (c_3)).
$$
Taking the constants sufficiently large, we may arrange that
this quantity is bounded by $1/2$. Then (\ref{eq:ineq_main}) trivially holds when 
$|u_1^{z_1}\cdots u_l^{z_l}u-1|\ge 1/2$, and without loss of generality we assume that
$|u_1^{z_1}\cdots u_l^{z_l}u-1|\le 1/2$.

Let $\log$ denote the principle value of the (complex) logarithm. 
There exists $z_0\in \Z$ such that $|z_0|\ll \|z\|$ and 
$$
T :=\log(u_1^{z_1}\cdots u_n^{z_n}u)=\pi i z_0+z_1\log u_1+\cdots+z_l\log u_l+\log u.
$$
It is convenient to set $u_0=-1$, so that $\log u_0=\pi i$.
Let $S:=u_1^{z_1}\cdots u_l^{z_l}u$. Since $|S-1|\le 1/2$,
$$
|T|=|\log S|\le 2|S-1|.
$$
Hence, it sufficient to establish a lower bound for $|T|$.
Note that since $S\ne 1$, we have $T\ne 0$.
For this purpose we use \cite[Cor.~10.1]{wald}, which we now recall.
We note that the result in \cite{wald} is stated using the logarithmic height
while here we use the exponential height. For simplicity, we take $E=e$ and $f=1$.

Let $D=[\Q(u_0,\ldots,u_l,u):\Q]$, $A_0,\ldots, A_l,B$ be numbers, greater than $e$, such that
\begin{align}\label{eq:ab}
&\hbox{H}(u_i)\le A_i,\;\; i=0,\ldots,l,\quad \hbox{H}(u)\le B, \nonumber\\
&\sum_{i=0}^l \frac{|\log u_i|}{\log A_i}+\frac{|\log u|}{\log
  B}\le e^{-1} (l+2)D.
\end{align}
We set
\begin{align*}
A&=\max\{A_0,\ldots, A_l,B\},\\
M&=\max_{i=0,\ldots,l}\left\{\frac{1}{\log A_i}+\frac{|z_i|}{\log B} \right\},\\
Z_0&=\max\left\{7+3\log(l+2), \log D\right\},\\
G_0&=\max\{4(l+2)Z_0,\log M, \log D\},\\
U_0 &= \max\left\{D^2\log A, D^{l+4} G_0Z_0(\log A_0)\cdots (\log A_l)(\log B)\right\}.
\end{align*}
Then according to \cite[Cor.~10.1]{wald},
\begin{equation}\label{eq:T}
|T|\ge \exp(-c\, U_0),
\end{equation}
where $c$ is an explicit positive constant depending only on $n$. We set $B:=c_2\hbox{H}(u)$ with
$c_2>1$. We note that
\[  |\log u|^2 \le \pi^2+\log ^2 |u|  \le \pi^2+[\Q(u):\Q]^2 \log ^2 \hbox{H}(u).\]
Therefore, taking $A_i$ sufficiently large, depending on $u_i$, and sufficiently large $c_2$,
we may arrange that (\ref{eq:ab}) holds. If $c_2$ is sufficiently large,  $A=B$.
Under the assumption (\ref{eq:cond}), we have $M\le \frac{c_3\|z\|}{\log B}$ with sufficiently large $c_3$
and also $G_0=\log M\le \log\left(\frac{c_3\|z\|}{\log B}\right)$.
Moreover, if $c_3$ is sufficiently large, then $U_0\ll \log\left(\frac{c_3\|z\|}{\log B}\right)\log B$.
Therefore, estimate (\ref{eq:T}) implies that
$$
|T|\ge \exp\left(-c_1\, \log\left(\frac{c_3\|z\|}{\log B}\right)\log B\right),
$$
where $c_1$ is an explicit positive constant. This completes the proof of the proposition.
\QED

\subsection{Equidistribution of box maps}  \label{Green-Tao dichotomy} 

A box map is an affine map 
$$
\iota:B:=[0,T_1]\times\cdots\times [0,T_k]\to \mathcal{L}(G)
$$
of the form
\begin{equation}\label{eq:box}
\iota: (t_1,\ldots,t_k)\mapsto v+t_1w_1+\cdots +t_kw_k
\end{equation}
with $v,w_1,\ldots,w_k\in \mathcal{L}(G)$.
We shall use the following result, which is a variation of our theorem
 \cite[Th.~2.1]{GS}, that implies equidistribution of box maps
under suitable Diophantine conditions. 
This result is based on the work of B.~Green and T.~Tao \cite{Green-Tao}.

We denote by $|B|$ the $k$-dimensional volume of the box $B$.

\begin{Theorem}\label{th:green-tao}
Fix $0<\theta\le 1$. There exist $L_1,L_2>0$ such that 
for every $\delta\in (0,\delta_0)$ and
every box map $\iota:B\to\mathcal{L}(G)$, 
 one of the following holds:
\begin{enumerate}
\item[(i)] For every $\theta$-\holder function $f:X\to \mathbb{R}$, $u\in \mathcal{L}(G)$, and $g\in G$, 
\begin{equation*}
\left|\frac{1}{|B|}\int_B f(\exp(u)\exp(\iota(t))g\Lambda)\, dt-\int_X f\, d\mu\right|\le \delta \|f\|_{ C^{\theta}}.
\end{equation*}
\item[(ii)] There exists $z\in \mathbb{Z}^d\backslash \{0\}$ such that
\begin{equation*}
\|z\|\ll \delta^{-L_1}\quad\hbox{and}\quad |\left<z,D\pi(w_i)\right>|\ll \delta^{-L_2}/T_i \quad \hbox{for
  all $i=1,\ldots,k$.}
\end{equation*}
\end{enumerate}
\end{Theorem}

\Proof 
In the case of Lipschitz functions $f$, this is \cite[Th.~2.1]{GS},
and the analogous result for H\"older functions can be deduced by a standard approximation
argument. Indeed, suppose that for some $f\in C^\theta(X)$, 
$u\in \mathcal{L}(G)$, and $g\in G$, 
\begin{equation}\label{eq:bad}
\left|\frac{1}{|B|}\int_B f(\exp(u)\exp(\iota(t))g\Lambda)\, dt-\int_X f\, d\mu\right|> \delta \|f\|_{ C^{\theta}}.
\end{equation}
Then one can find a Lipschitz function $f_\varepsilon$ such that
$$
 \|f_{\varepsilon } -f \|_{C^0} \leq \varepsilon ^{\theta} \|f\|_{C^\theta}  \quad \text{and} \quad
   \|f_\varepsilon\|_{Lip} \ll \varepsilon^{-\dim(X)-1} \|f\|_{C^0}
$$
(see, for instance, \cite[Lem.~2.4]{GS}). Then 
taking $\varepsilon=(\delta/3)^{1/\theta}$,
we deduce from (\ref{eq:bad}) that
\begin{align*}
\left|\frac{1}{|B|}\int_B f_\varepsilon(\exp(u)\exp(\iota(t))g\Lambda)\, dt-\int_X f_\varepsilon\, d\mu\right|&> 
(\delta-2\varepsilon^\theta) \|f\|_{ C^{\theta}}\\
&\gg 
\varepsilon^{\dim(X)+1} (\delta-2\varepsilon^\theta) \|f_\varepsilon\|_{Lip}\\
&\gg \delta^{(\dim(X)+1)/\theta+1} \|f_\varepsilon\|_{Lip}.
\end{align*}
Now the theorem for Lipschitz functions implies that (ii) holds with some $L_1,L_2>0$ depending on $\theta$.
\QED

\subsection{3-mixing under irreducibility condition}\label{sec:3_irred}

The action of $D\alpha(\Z^l)$ preserves the rational structure on $\mathcal{L}(G)$
defined by the lattice $\Lambda$. In particular, it follows that each character
$\chi$ in $\mathcal{X}(\alpha)$ is of the form $\chi(z)=u_1^{z_1}\cdots u_l^{z_l}$
where $u_i$'s are algebraic numbers. The Galois group $\hbox{Gal}(\bar \Q/\Q)$
naturally acts on $\mathcal{X}(\alpha)$ and on $\mathcal{L}(G)_{\bar\Q}$.
We fix an orbit $\mathcal{X}_0\subset \mathcal{X}'(\alpha)$ of the Galois group
and for each $\chi\in \mathcal{X}_0$, we fix a vector $w_\chi\in \mathcal{L}_\chi$
whose coordinates are algebraic integers, so that the vectors
$w_\chi$ are also conjugate under the action of the Galois group.
Let $W_{\C}$ be the Lie subalgebra of $\mathcal{L}(G)\otimes \C$
generated by the vectors $w_\chi$, $\chi\in \mathcal{X}_0$,
and $W=W_{\C}\cap \mathcal{L}(G)$. We also fix a basis $\{w_i\}$ of $W$.

In this section, we prove Theorem \ref{th:3_mixing} under the irreducibility assumption: $D\pi(W)$ is not
contained in a proper rational subspace. Let $\bar w_\chi=D\pi(w_\chi)$, $\chi\in \mathcal{X}_0$.
We observe that under this assumption the coordinates of each the vectors $\bar w_\chi$
are linearly independent over $\Q$. Indeed, if we suppose that $\left<a,\bar w_\chi\right>=0$
for some $a\in \Q^d\backslash \{0\}$, then applying the action of the Galois group, 
we deduce that $\left<a,\bar w_\chi\right>=0$ for all $\chi\in \mathcal{X}_0$.
Since $D\pi(W)$ is spanned over $\C$ by the vectors $w_\chi$, $\chi\in \mathcal{X}_0$,
this would imply that $D\pi(W)$ is contained in a proper rational subspace, which contradicts
our assumption.

Let 
\begin{equation}\label{eq:N}
N=N(z_1,z_2,z_3):=\exp\left(\min_{i\ne j} \|z_i-z_j\|\right).
\end{equation}
Without loss of generality,
we may assume that $z_0=0$ and $N=\exp(\|z_1-z_2\|)$.
We set $\epsilon=N^{-\kappa}$, where $\kappa>0$ is a fixed parameter
which is sufficiently small and will be specified later
(see (\ref{eq:kappa1}), (\ref{eq:kappa2}), (\ref{eq:kappa3}), (\ref{eq:kappa4}) below).

We fix a fundamental domain $F\subset G$ for $X=G/\Lambda$ and set $E=\exp^{-1}(F)$.
As in \cite[Section 3]{GS}, we may arrange that $E$ is bounded and has piecewise smooth boundary.
Since the Haar measure on $G$ is the image under $\exp$ of a suitably normalised
Lebesgue measure on $\mathcal{L}(G)$ \cite[1.2.10]{CG}, we obtain
\begin{align}\label{eq:000}
&\int_X f_0(x)f_1(\alpha(z_1)x)f_2(\alpha(z_2)x)\,d\mu(x)\\
=&\int_E f_0(\exp(u)\Lambda)f_1(\exp(D\alpha(z_1)u)\Lambda)f_2(\exp(D\alpha(z_2)u)\Lambda)\, du.\nonumber
\end{align}
We choose a basis of $\mathcal{L}(G)$ that contains the basis $\{w_i\}$ of $W$ 
and tessellate $\mathcal{L}(G)$ by cubes $C$ of size $\epsilon$ with respect to this basis.
Since $E$ has piecewise smooth boundary, we obtain 
\begin{equation}\label{eq:vol}
\left|E-\bigcup_{C\subset E} C\right|\ll \epsilon,
\end{equation}
and
\begin{align}\label{eq:ee1}
&\int_E f_0(\exp(u)\Lambda)f_1(\exp(D\alpha(z_1)u)\Lambda)f_2(\exp(D\alpha(z_2)u)\Lambda)\, du\\
=&\sum_{C\subset E} 
\int_C f_0(\exp(u)\Lambda)f_1(\exp(D\alpha(z_1)u)\Lambda)f_2(\exp(D\alpha(z_2)u)\Lambda)\, du\nonumber\\
&+O(\epsilon\|f_0\|_{C^\theta}\|f_1\|_{C^\theta}\|f_2\|_{C^\theta}).\nonumber
\end{align}
For every cube $C$, we pick a point $u_C\in C$. Then since $f_0$ is $\theta$-H\"older,
\begin{align}\label{eq:ee2}
&\int_C f_0(\exp(u)\Lambda)f_1(\exp(D\alpha(z_1)u)\Lambda)f_2(\exp(D\alpha(z_2)u)\Lambda)\, du\\
=& f_0(\exp(u_C)\Lambda) \int_C f_1(\exp(D\alpha(z_1)u)\Lambda)f_2(\exp(D\alpha(z_2)u)\Lambda)\, du\nonumber\\
&+O(\epsilon^\theta \|f_0\|_{C^\theta}\|f_1\|_{C^\theta}\|f_2\|_{C^\theta}).\nonumber
\end{align}
We decompose each cube $C$ as $C=B'+B$, where $B$ is a cube in $W$ and $B'$
is a cube in the complementary subspace.  

We claim that for sufficiently small $\kappa>0$ and all sufficiently large $N$
defined in (\ref{eq:N}),
\begin{align}\label{eq:claim}
&\frac{1}{|B|} \int_B f_1(\exp(v_1+D\alpha(z_1)b)\Lambda)
f_2(\exp(v_2+D\alpha(z_2)b)\Lambda)\, db\\
 =&
\left( \int_X f_1\, d\mu\right)\left( \int_X f_2\, d\mu\right)
+O(N^{-\kappa}\|f_1\|_{C^\theta} \|f_2\|_{C^\theta}),\nonumber
\end{align}
uniformly over the cubes $B$ and $v_1,v_2\in \mathcal{L}(G)$,

Suppose first that (\ref{eq:claim}) holds. Then using uniformity over $v_1,v_2$, we deduce that
\begin{align*}
&\frac{1}{|C|}\int_C f_1(\exp(D\alpha(z_1)u)\Lambda)f_2(\exp(D\alpha(z_2)u)\Lambda)\,du\\
=&\frac{1}{|B'||B|}\int_{B'}\int_B f_1(\exp(D\alpha(z_1)b'+D\alpha(z_1)b)\Lambda)f_2(\exp(D\alpha(z_2)b'+D\alpha(z_2)b)\Lambda)\,dbdb'\\
=&\left(\int_X f_1\, d\mu\right)\left(\int_X f_2\, d\mu\right)+O(N^{-\kappa}\|f_1\|_{C^\theta}\|f_2\|_{C^\theta}).
\end{align*}
Combining this estimate with (\ref{eq:ee1})--(\ref{eq:ee2}), we obtain
\begin{align*}
&\int_E f_0(\exp(u)\Lambda)f_1(\exp(D\alpha(z_1)u)\Lambda)f_2(\exp(D\alpha(z_2)u)\Lambda)\, du\\
=&\left(\sum_{C\subset E} f_0(\exp(u_C)\Lambda)|C| \right)\left(\int_X f_1\, d\mu\right)
\left(\int_X f_2\, d\mu\right)\\
&+O\left((N^{-\kappa}+\epsilon^\theta)\|f_0\|_{C^\theta}\|f_1\|_{C^\theta}\|f_2\|_{C^\theta}\right).
\end{align*}
Since $f$ is $\theta$-H\"older and (\ref{eq:vol}) holds,
\begin{align*}
&\sum_{C\subset E} f_0(\exp(u_C)\Lambda)|C|=\sum_{C\subset E} \int_C
f_0(\exp(u)\Lambda)\, du +O\left(\epsilon^\theta\|f_0\|_{C^\theta}\right)\\
=&\int_E
f_0(\exp(u)\Lambda)\, du+O\left((\epsilon+\epsilon^\theta)\|f_0\|_{C^\theta}\right)
=\int_X
f_0\,d\mu +O\left(\epsilon^\theta\|f_0\|_{C^\theta}\right).
\end{align*}
Hence, 
\begin{align}\label{eq:999}
&\int_E f_0(\exp(u)\Lambda)f_1(\exp(D\alpha(z_1)u)\Lambda)f_2(\exp(D\alpha(z_2)u)\Lambda)\, du\\
=& \left(\int_X f_0\, d\mu\right) \left(\int_X f_1\, d\mu\right)
\left(\int_X f_2\, d\mu\right)
+O\left(N^{-\kappa\theta}\|f_0\|_{C^\theta}\|f_1\|_{C^\theta}\|f_2\|_{C^\theta}\right).\nonumber
\end{align}
This proves the required estimate when $N$ is sufficiently large, and it is also clear that this
estimate holds for $N$ in bounded intervals. Hence, Theorem \ref{th:3_mixing} follows.
Now it remains to prove the claim (\ref{eq:claim}).

To prove (\ref{eq:claim}), we apply 
Theorem \ref{th:green-tao}   to the nilmanifold $X\times X=(G\times G)/(\Lambda\times
\Lambda)$ with $\delta=N^{-\kappa}$. We assume that $N$
is sufficiently large, so that Theorem \ref{th:green-tao} applies.
Let $f=f_1\otimes f_2$. Clearly,
$$
\int_{X\times X} f\,d(\mu\otimes\mu)=\left(\int_X f_1\, d\mu\right)
\left(\int_X f_1\, d\mu\right) \quad\hbox{and}\quad
\|f\|_{C^\theta}\ll \|f_1\|_{C^\theta}\|f_2\|_{C^\theta}.
$$
We consider the map $\iota:[0,\epsilon]^k\to \mathcal{L}(G)$ defined by
$$
\iota(t)=\left(v'_1+\sum_{i=1}^k t_i D\alpha(z_1) w_i,v'_2+\sum_{i=1}^k t_i D\alpha(z_2) w_i\right)
$$
with suitably chosen  $v_1',v_2'\in \mathcal{L}(G)$, so that
$$
\int_B f_1(\exp(v_1+D\alpha(z_1)b)\Lambda)
f_2(\exp(v_2+D\alpha(z_2)b)\Lambda)\, db
=\int_{[0,\epsilon]^k} f(\iota(t)\Lambda)\, dt,
$$
It is sufficient to prove that 
$$
\epsilon^{-k}\int_{[0,\epsilon]^k} f(\iota(t)\Lambda)\, dt=
\int_{X\times X} f\,d(\mu\otimes\mu)+O(\delta\|f\|_{C^\theta}).
$$
Applying 
Theorem \ref{th:green-tao}, we deduce that either
\begin{equation}\label{eq:b0}
\left|\epsilon^{-k}\int_{[0,\epsilon]^k} f(\iota(t)\Lambda)\, dt-\int_{X\times X} f\,d(\mu\otimes\mu)\right|\le \delta\|f\|_{C^\theta},
\end{equation}
or there exists $(a_1,a_2)\in (\Z^d)^2\backslash \{(0,0)\}$ such that
\begin{equation}\label{eq:b1}
\max(\|a_1\|,\|a_2\|)\ll \delta^{-L_1}=N^{\kappa L_1}
\end{equation}
and 
\begin{equation}\label{eq:b2}
\left| \left<a_1,(D\pi)D\alpha(z_1)w_i\right>
+\left<a_2,(D\pi)D\alpha(z_2)w_i\right>\right|\ll \delta^{-L_2}/\epsilon=N^{\kappa (L_2+1)}
\end{equation}
for all $i=1,\ldots,k$. 

We shall show that if $\kappa>0$ is sufficiently small and $N$ is sufficiently large,
then (\ref{eq:b1})--(\ref{eq:b2}) fails. Suppose that (\ref{eq:b1})--(\ref{eq:b2}) holds.
Since each of the vectors $w_\chi$, $\chi\in\mathcal{X}_0$, is a linear combination of vectors
$w_i$, we deduce from (\ref{eq:b2}) that
\begin{equation}\label{eq:b220}
\left| \left<a_1,(D\pi)D\alpha(z_1)w_\chi\right>
+\left<a_2,(D\pi)D\alpha(z_2)w_\chi\right>\right|\ll N^{\kappa (L_2+1)}
\quad\hbox{for all $\chi\in\mathcal{X}_0$.}
\end{equation}
Since $D\alpha(z)w_\chi=\chi(z)w_\chi$ and $\bar w_\chi= D\pi(w_\chi)$,
(\ref{eq:b220}) becomes
\begin{equation}\label{eq:b22}
\left| \chi(z_1)\left<a_1, \bar w_\chi\right>
+\chi(z_2)\left<a_2,\bar w_\chi\right>\right|\ll N^{\kappa (L_2+1)}\quad\hbox{for all $\chi\in\mathcal{X}_0$.} 
\end{equation}
We divide the argument into three cases.

{\it Case 1:} suppose that $a_1=0$. Then $a_2\ne 0$ and $\left<a_2,\bar w_\chi \right>\ne 0$.
Moreover, by \cite[Th.~7.3.2]{BG}, 
\begin{equation}\label{eq:n1}
|\left<a_2,\bar w_\chi \right>|\gg \|a_2\|^{-d-1}\gg N^{-\kappa L_1(d+1)}.
\end{equation}
By Lemma \ref{l:expansion}, there exists $\chi\in\mathcal{X}_0$ such that $|\chi(z_2)|\ge N^c$
with fixed $c>0$. Hence, it follows from (\ref{eq:b22}) that
\begin{equation}\label{eq:n2}
|\left<a_2,\bar w_\chi \right>|\ll N^{\kappa(L_2+1)-c}.
\end{equation}
We assume that the parameter $\kappa>0$ satisfies
\begin{equation}\label{eq:kappa1}
-\kappa L_1(d+1)>\kappa(L_2+1)-c.
\end{equation}
Comparing (\ref{eq:n1}) and (\ref{eq:n2}), we get a contradiction if $N$
is sufficiently large. Hence, we may assume that $a_1\ne 0$.

{\it Case 2:} suppose that $a_1\ne 0$ and $\chi(z_1)\left<a_1, \bar w_\chi\right>+\chi(z_2)\left<a_2,\bar
  w_\chi\right>= 0$ for some $\chi\in \mathcal{X}_0$. 
Since the Galois group acts transitively on the set $\mathcal{X}_0$,
it follows that this equality holds for all $\chi\in \mathcal{X}_0$.
By Lemma \ref{l:expansion}, there exists $\chi\in \mathcal{X}_0$ such that $|\chi(z_2-z_1)|\ge N^c$
with fixed $c>0$. Then
\begin{equation}\label{eq:contr}
|\left<a_2, \bar w_\chi\right>|=|\chi(z_1-z_2)| |\left<a_1,\bar
  w_\chi\right>|\ge N^c |\left<a_1,\bar w_\chi\right>|.
\end{equation}
Since $a_1\ne 0$, we have $\left<a_1, \bar w_\chi\right>\ne 0$, and by \cite[Th.~7.3.2]{BG},
$$
|\left<a_1, \bar w_\chi\right>|\gg \|a_1\|^{-d-1}\gg N^{-\kappa L_1(d+1)}.
$$
On the other hand,
$$
|\left<a_2, \bar w_\chi\right>|\ll \|a_2\|\ll N^{\kappa L_1}.
$$
Hence, we deduce that 
$$
N^{-\kappa L_1(d+1)+c}\ll N^{\kappa L_1}.
$$
We choose the parameter $\kappa>0$ so that
\begin{equation}\label{eq:kappa2}
-\kappa L_1(d+1)+c> \kappa L_1.
\end{equation}
Then when $N$ is sufficiently large, we get a contradiction.

{\it Case 3:} suppose that $a_1\ne 0$ and $\chi(z_1)\left<a_1, \bar w_\chi\right>+\chi(z_2)\left<a_2,\bar
  w_\chi\right>\ne 0$ for all $\chi\in\mathcal{X}_0$. This is the most difficult part of the proof.

Since $a_1\ne 0$, we have $\left<a_1, \bar w_\chi\right>\ne 0$, and by \cite[Th.~7.3.2]{BG},
$$
|\left<a_1, \bar w_\chi\right>|\gg \|a_1\|^{-d-1}\gg N^{-\kappa L_1(d+1)}. 
$$
We set $u=-\frac{\left<a_2, \bar w_\chi\right>}{\left<a_1, \bar w_\chi\right>}$. 
By Lemma \ref{l:expansion},
there exists $\chi\in \mathcal{X}_0$ such that $|\chi(z_1)|\gg N^{c}$
with fixed $c>0$.
It follows from (\ref{eq:b22}) that for this $\chi$, we have the estimate
\begin{equation}\label{eq:e1}
|\chi(z_2-z_1)u-1|\ll \frac{N^{\kappa (L_2+1)}}{|\chi(z_1)| |\left<a_1, \bar w_\chi\right>|}
\ll N^{\kappa (L_2+1+ L_1(d+1))-c}.
\end{equation}
Let $K_1:=\kappa (L_2+1+ L_1(d+1))-c$.

Next, we compare this estimate with the lower estimate provided by  Proposition \ref{p:waldschmidt}.
We note that
\begin{align*}
\hbox{H}(u)=\prod_v \max (|\left<a_1, \bar w_\chi\right>|_v,|\left<a_2, \bar w_\chi\right>|_v)^{1/[\Q(u):\Q]}.
\end{align*} 
For all non-Archimedian places $v$,
$$  
|\left<a_i, \bar w_\chi\right>|_v\le 1,
$$
and for all Archimedian $v$,
$$
|\left<a_i, \bar w_\chi\right>|_v\ll \|a_i\|\ll N^{\kappa L_1}.
$$ 
Therefore,
\begin{align}\label{eq:hh}
\hbox{H}(u) \ll  N^{K_2},
\end{align}
where $K_2:=\kappa L'_1$ with fixed $L_1'>0$.
We take the parameter $\kappa>0$ so that
\begin{equation}\label{eq:kappa3}
K_2=\kappa L_1'<1.
\end{equation}
Then assuming that $N$ is sufficiently large, we obtain
\begin{equation}\label{eq:ccc}
\log (c_2 \hbox{H}(u))\le \log (c'_2 N^{K_2})\le \log N,
\end{equation}
where $c_2'>1$ depends on the implicit constant in the estimate (\ref{eq:hh}).
We recall that we have chosen the indices, so that
$$
\log N=\|z_2-z_1\|.
$$
Since (\ref{eq:ccc}) holds, Proposition \ref{p:waldschmidt} applies, and we deduce that 
\begin{align*}
|\chi(z_2-z_1)u-1| &\ge 
\exp\left(-
  c_1\log(c_2\hbox{H}(u))\log\left(\frac{c_3\|z_2-z_1\|}{\log(c_2\hbox{H}(u))}\right)\right).
\end{align*}
Without loss of generality, we may assume that $c_3>e$. 
Since the function $x\mapsto x\log(C/x)$ is increasing for $x\le C/e$, we deduce that
\begin{align}\label{eq:e2}
|\chi(z_2-z_1)u-1| &\ge  \exp\left(-c_1\log(c'_2N^{K_2})\log\left(\frac{c_3\log
      N}{\log(c_2'N^{K_2})}\right)\right) \\
&\ge  \exp\left(-c_1\log(c'_2N^{K_2})\log\left(c_3K_2^{-1}\right)\right).\nonumber
\end{align} 
Comparing (\ref{eq:e1}) and (\ref{eq:e2}), we conclude that
\begin{equation}\label{eq:est0}
K'_2\log N+M_2\le K_1\log N+M_1,
\end{equation}
where $K_2':=-c_1K_2\log(c_3 K_2^{-1})$,
$M_2 : =-c_1\log(c_2')\log(c_3K_2^{-1})$, and $M_1$ is determined 
by the implicit constant in (\ref{eq:e1}).
We observe that as $\kappa\to 0^+$, we have $K_2'\to 0^-$ and $K_1\to -c<0$.
Therefore, taking the parameter $\kappa>0$ sufficiently small, we may arrange that
\begin{equation}\label{eq:kappa4}
K_2'>K_1.
\end{equation}
Then when $N$ is sufficiently large, (\ref{eq:est0}) fails.
This shows that (\ref{eq:b1})--(\ref{eq:b2}) fails and (\ref{eq:b0}) holds
when $N$ is sufficiently large. Now we have verified the claim (\ref{eq:claim})
and completed the proof of Theorem \ref{th:3_mixing} under the irreducibility condition.

In order to prove Theorem \ref{th:3_mixing} in general, we observe that using the same argument,
one can deduce the following more general version of the estimate (\ref{eq:claim}):
for all sufficiently large $N$ defined in (\ref{eq:N}),
\begin{align}\label{eq:claim_2}
&\frac{1}{|B|} \int_B f_1(h_1\beta_1(\exp(v_1+D\alpha(z_1)b))\Lambda)
f_2(h_2\beta_2(\exp(v_2+D\alpha(z_2)b))\Lambda)\, db\\
 =&
\left( \int_X f_1\, d\mu\right)\left( \int_X f_2\, d\mu\right)
+O(N^{-\kappa}\|f_1\|_{C^\theta} \|f_2\|_{C^\theta})\nonumber
\end{align}
uniformly over the cubes $B$, $h_1,h_2\in G$, $v_1,v_2\in \mathcal{L}(G)$,
and automorphisms $\beta_1,\beta_2$ of $G$ which act trivially on $G/G'$.
Indeed, 
\begin{align*}
&\int_B f_1(h_1\beta_1(\exp(v_1+D\alpha(z_1)b))\Lambda)f_2(h_2\beta_2(\exp(v_2+D\alpha(z_2)b))\Lambda)\,
db\\
=&
\int_B f_1(h_1\exp(D\beta_1(v_1)+D\beta_1D\alpha(z_1)b)\Lambda)
f_2(h_2\exp(D\beta_2(v_2)+D\beta_2D\alpha(z_2)b))\Lambda)\, db,
\end{align*}
and to prove (\ref{eq:claim_2}), we can apply Theorem \ref{th:green-tao}
to the map
$$
\iota:t\mapsto \left(v'_1+\sum_{i=1}^k t_i D\beta_1D\alpha(z_1)w_i, v_2'+\sum_{i=1}^k t_i D\beta_2 D\alpha(z_2)w_i\right)
$$
As in the above proof, either (\ref{eq:claim_2}) holds, or an analogue of
(\ref{eq:b2}) holds, but since $D\pi D\beta_i=D\pi$, this reduces to the exactly same estimate as
(\ref{eq:b2}). Therefore, (\ref{eq:claim_2}) follows. Now we combine (\ref{eq:claim_2})
with the argument (\ref{eq:000})--(\ref{eq:999}) to deduce that 
\begin{align}\label{eq:final}
&\int_X f_0(x)f_1(h_1\beta_1(\alpha(z_1)(x)))f_2(h_2\beta_2(\alpha(z_2)(x)))\, d\mu(x)\\
=& \left(\int_X f_0\, d\mu\right) \left(\int_X f_1\, d\mu\right)
\left(\int_X f_2\, d\mu\right)
+O\left(N^{-\kappa\theta}\|f_0\|_{C^\theta}\|f_1\|_{C^\theta}\|f_2\|_{C^\theta}\right)\nonumber
\end{align}
uniformly over $h_1,h_2\in G$ and automorphisms $\beta_1,\beta_2$ of $G$ that preserve $\Lambda$
act trivially on $G/G'$. We will use this estimate to establish Theorem \ref{th:3_mixing} in general.

\subsection{3-mixing in general}\label{sec:3_gen}
Let $W$ be the Lie subalgebra of $\mathcal{L}(G)$ introduced in Section \ref{sec:3_irred}.
By \cite[Ch.~5, Sec.~5]{Starkov}, there exists a closed connected normal subgroup $M$ of $G$ such that
$M/(M\cap \Lambda)$ is compact, and 
$$
\overline{\exp(W)g\Lambda}=Mg\Lambda\quad\hbox{for almost every $g\in G$.}
$$
Since we may replace the lattice $\Lambda$ by its conjugate, we assume that 
$$
\overline{\exp(W)\Lambda}=M\Lambda.
$$
We note that the group $M$ satisfies the following properties:
\begin{enumerate}
\item[(i)] $M$ is $\alpha(\Z^l)$-invariant,
\item[(ii)] $D\pi(W)$ is not contained in a proper rational subspace of $\mathcal{L}(M/M')$,
\item[(iii)] $[G,M]\subset M'$.
\end{enumerate}
Properties (i)--(iii) can be verified exactly as in the proof of \cite[Lem.~3.4]{GS}.

We give the proof of Theorem \ref{th:3_mixing} using induction on dimension of $X$.
For this, we use that $X=G/\Lambda$ fibers over $Y=G/(M\Lambda)$ with fibers isomorphic
to 
$$
R=M\Lambda/\Lambda\simeq M/(M\cap\Lambda).
$$
The invariant measure on $X$ can be decomposed as
$$
\int_X f\, d\mu=\int_Y \int_R f(yr)\, d\mu_R(r) d\mu_Y(y),\quad f\in C(X),
$$
where $\mu_Y$ and $\mu_R$ are normalised invariant measure on $Y$ and $R$ respectively.
Since the fibration is $\alpha(\Z^l)$-equivariant (by (i)),
\begin{align}\label{eq:induction0}
&\int_X f_0(x)f_1(\alpha(z_1)x)f_2(\alpha(z_2)x)\, d\mu(x)\\
=&\int_Y \left(\int_R  f_0(yr)f_1(\alpha(z_1)(y) \alpha(z_1)(r))f_2(\alpha(z_2)(y) \alpha(z_2)(r))  \,d\mu_R(r)\right)d\mu_Y(y)\nonumber\\
=&\int_F \left(\int_R  f_0(gr)f_1(\alpha(z_1)(g) \alpha(z_1)(r))f_2(\alpha(z_2)(g) \alpha(z_2)(r))  \,d\mu_R(r)\right)dm_F(g),\nonumber
\end{align}
where $F\subset G$ is a bounded fundamental domain for $G/(M\Lambda)$, and $m_F$ is the measure on $F$
induced by $\mu_Y$. We shall show that for $N$ defined in (\ref{eq:N})
and some $\eta>0$,
\begin{align}\label{eq:induction}
&\int_R  f_0(gr)f_1(\alpha(z_1)(g) 
\alpha(z_1)(r)) f_2(\alpha(z_2)(g) \alpha(z_2)(r))  \, d \mu_R (r)\\
=& \left(\int_R f_0(gr)d\mu_R(r)\right) \left(\int_R f_1(\alpha(z_1)(g)r)d\mu_R(r)\right)
\left(\int_R
  f_2(\alpha(z_2)(g)r)d\mu_R(r)\right)\nonumber\\
&+O(N^{-\eta}\|f_0\|_{C^\theta}\|f_1\|_{C^\theta}\|f_2\|_{C^\theta}) \nonumber
\end{align}
uniformly over $g\in F$.

Suppose that (\ref{eq:induction}) holds. Then combining (\ref{eq:induction0}) and (\ref{eq:induction}),
we obtain
\begin{align*}
&\int_X f_0(x)f_1(\alpha(z_1)x)f_2(\alpha(z_2)x)\, d\mu(x)\\
=&\int_Y \bar f_0(y)\bar f_1(\alpha(z_1)y)\bar f_2(\alpha(z_2)y)\, d\mu_Y(y)
+O(N^{-\eta}\|f_0\|_{C^\theta}\|f_1\|_{C^\theta}\|f_2\|_{C^\theta}),
\end{align*}
where the functions $\bar f_i$ on $Y$ are defined by 
$$
y\mapsto \int_R f_i(yr)\,d\mu_R(r).
$$
Since $\dim(Y)<\dim(X)$, it follows from the inductive assumption that for some $\eta>0$,
\begin{align*}
&\int_Y \bar f_0(y)\bar f_1(\alpha(z_1)y)\bar f_2(\alpha(z_2)y)\, d\mu_Y(y)\\
=& \left(\int_Y \bar f_0\, d\mu_Y\right) 
\left(\int_Y \bar f_1\, d\mu_Y\right) \left(\int_Y \bar f_2\, d\mu_Y\right) 
+O(N^{-\eta}\|\bar f_0\|_{C^\theta}\|\bar f_1\|_{C^\theta}\|\bar f_2\|_{C^\theta})\\
=& \left(\int_X f_0\, d\mu\right) 
\left(\int_X f_1\, d\mu\right) \left(\int_X f_2\, d\mu\right) 
+O(N^{-\eta}\|f_0\|_{C^\theta}\|f_1\|_{C^\theta}\|f_2\|_{C^\theta}),
\end{align*}
and this completes the proof of Theorem \ref{th:3_mixing}. Hence, it remains to prove (\ref{eq:induction}).

To prove (\ref{eq:induction}), we write
$$
\alpha(z_i)(g)=a_im_i\lambda_i\quad\hbox{with $a_i\in F$, $m_i\in M$, $\lambda_i\in\Lambda$,}\quad
 i=1,2.
$$
Then 
\begin{align*}
&\int_R  f_0(gr)f_1(\alpha(z_1)(g) 
\alpha(z_1)(r)) f_2(\alpha(z_2)(g) \alpha(z_2)(r))  \, d \mu_R (r)\\
=&\int_R  f_0(gr)f_1(a_1m_1 
\beta_1(\alpha(z_1)(r))) f_2(a_2m_2\beta_2(\alpha(z_2)(r)))  \, d \mu_R (r),
\end{align*}
where $\beta_i$'s are the maps induced by the automorphisms $m\mapsto \lambda_i m\lambda_i^{-1}$.
We observe that because of (ii), $W\subset \mathcal{L}(M)$ satisfies the irreducibility assumption
of Section \ref{sec:3_irred}, and by (iii), the automorphisms $\beta_i$ act trivially on $M/M'$.
Hence, (\ref{eq:final}) holds. We apply (\ref{eq:final}) to the functions on $R$ defined by
$$
\phi_0(r):=f_0(gr)\quad\hbox{and}\quad \phi_i(r):=f_i(a_ir),\quad\quad i=1,2.
$$
This gives
\begin{align*}
&\int_R \phi_0(r)\phi_1(m_1\beta_1(\alpha(z_1)r))\phi_2(m_2\beta_2(\alpha(z_2)r))\, d\mu_R(r)\\
=& \left(\int_R\phi_0\, d\mu_R\right)\left(\int_R\phi_1\, d\mu_R\right)\left(\int_R\phi_2\, d\mu_R\right)
+O(N^{-\eta}\|\phi_0\|_{C^\theta}\|\phi_1\|_{C^\theta}\|\phi_2\|_{C^\theta})\\
=& \left(\int_Rf_0(gr)\, d\mu_R(r)\right)\left(\int_Rf_1(\alpha(z_1)(g)r)\, d\mu_R(r)\right)
\left(\int_R f_2(\alpha(z_2)(g)r)\, d\mu_R(r)\right)\\
&+O(N^{-\eta}\|f_0\|_{C^\theta}\|f_1\|_{C^\theta}\|f_2\|_{C^\theta}).
\end{align*}
This implies (\ref{eq:induction}) and completes the proof of Theorem \ref{th:3_mixing}.


\section{Higher order mixing}

 The aim of this section is to prove Theorem \ref{th:multmixing}.
We shall use notation introduced in Section \ref{sec:notation}.
In Section \ref{sec:dioph2} we prepare Diophantine estimates.
Then in Section \ref{sec:mult_irr} we give a proof of Theorem \ref{th:multmixing} under an irreducibility
condition, and in Section \ref{sec:mult_gen} we give a proof in general using an inductive argument.

We note that it is sufficient to prove Theorem \ref{th:multmixing}
for a collection of functions $f_i\in L^\infty(X)$ which is dense in $L^1(X)$. 
Hence, we may assume that $f_0,\ldots,f_s\in C^\theta(X)$.
Furthermore, we may assume that $z_0=0$.

\subsection{Diophantine estimates}\label{sec:dioph2}

Let $K$ be a number field and $S$ a finite set of places of $K$ containing all the archimedean places.
We denote by $U_S$ the ring of $S$-units, namely, the groups of elements $x$ in $K$ such that
$|x|_v=1$ for $v\notin S$.  For a vector $\bar x\in K^s$, we define its (relative) height by
$$
\hbox{H}(\bar x)=\prod_{v} \max(1,\|\bar x\|_v),
$$
where $v$ runs  the set of all places of $K$, and $\|\bar x\|_v=\max_i |x_i|_v$.

\begin{proposition}\label{p:subspace}
Let $v\in S$ and $b_1,\ldots b_s\in K\backslash \{0\}$. Then for every $\epsilon>0$, the inequality
\begin{equation}\label{eq:subspace}
\left|b_1 +\sum_{j=2}^s b_{j}  x_j\right|_v<\hbox{H}(\bar x)^{-\epsilon},
\end{equation}
has finitely many solutions $\bar x\in U_S$ such that no proper subsum of 
$b_1+\sum_{j=2}^s b_j x_j$ vanishes.
\end{proposition}

We call such solutions of (\ref{eq:subspace}) {\em nondegenerate}.

\vspace{0.2cm}

\proof
We prove the proposition by induction on $s$. Note that when $s=1$, the statement
holds trivially because there are only finitely many solutions of $\hbox{H}(\bar x)<c$.

Given a solution $\bar x$ of (\ref{eq:subspace}), we set $\bar y=(1,\bar x)$, and
we denote by $j_v=j_v(\bar x)$ the first index $j_v$  such that 
$$
|y_{j_v}|_v=\|\bar y\|_v\ge 1.
$$
Partitioning the set of solutions according to the index $j_v$,
we may assume that this index is fixed.

We introduce a family of linear forms $L_{wj}(\bar y)$, with $w\in S$ and $j=1,\ldots, s$,  defined by
\begin{align*}
&L_{wj}(\bar y)=y_j,\quad\quad (w,j)\ne (v,j_v),\\
&L_{wj}(\bar y)=\sum_{j=1}^s b_j y_j,\quad (w,j)= (v,j_v).
\end{align*}
Then if $\bar y=(1,\bar x)$ corresponds to a solution of (\ref{eq:subspace}), 
\begin{align*}
&\prod_{j=1}^s |L_{wj}(\bar y)|_w= \prod_{j=1}^s |y_j|_w,\quad w\ne v,\\
&\prod_{j=1}^s |L_{vj}(\bar y)|_v= |L_{v j_v}(\bar y)|_v \prod_{j\ne j_v} |y_j|_v
< \hbox{H}(\bar y)^{-\epsilon} \prod_{j=1}^s |y_j|_v,
\end{align*}
and by the product formula,
\begin{align}\label{eq:subspace2}
\prod_{w\in S}\prod_{j=1}^s |L_{wj}(\bar y)|_w < \hbox{H}(\bar y)^{-\epsilon}.
\end{align}
By the W.~Schmidt subspace theorem \cite[Cor.~7.2.5]{BG}, all the solutions of 
(\ref{eq:subspace2}) are contained in a finite union of proper linear subspaces of $K^s$.
Partitioning solutions of (\ref{eq:subspace}) according to theses subspaces,
we may assume that these solutions additionally satisfy a nontrivial linear relation
\begin{equation}\label{eq:lin}
c_1+\sum_{j=2}^s c_j x_j=0
\end{equation}
with $c_1,\ldots,c_s\in K$. 

Suppose that $c_1\ne 0$. Given a solution $\bar x$ of (\ref{eq:lin}), we pick a minimal
$J\subset\{2,\ldots,s\}$ such that 
\begin{equation}\label{eq:lin1}
c_1+\sum_{j\in J} c_j x_j=0.
\end{equation}
Then no proper subsum in (\ref{eq:lin1}) vanishes. 
It follows from the finiteness of the number of non-degenerate solutions of unit equations
\cite[Th.~7.4.2]{BG} that $x_{j}$, $j\in J$, varies over a finite set.
This shows that for every solution $\bar x$ of (\ref{eq:lin})
there exists $j_0=2,\ldots, s$ such that $x_{j_0}$ belongs to a fixed finite set.
Hence, in order to prove finiteness of nondegenerate solutions (\ref{eq:subspace}), we may assume that
$x_{j_0}$ is fixed. Then (\ref{eq:subspace}) becomes
\begin{equation}\label{eq:subspace3}
\left|(b_1+b_{j_0} x_{j_0}) +\sum_{j\ne j_0} b_j x_j\right|_v<\hbox{H}(\bar x)^{-\epsilon}.
\end{equation}
Since we are assuming that no proper subsum in (\ref{eq:subspace}) vanishes,
$b_1+b_{j_0} x_{j_0}\ne 0$ and no proper subsum in (\ref{eq:subspace3}) vanishes as well.
Let $\bar x'=(x_j:j\ne j_0)$. Then 
$\hbox{H}(\bar x')\le \hbox{H}(\bar x)$.
Hence, by the inductive assumption,
the number of nondegenerate solutions $\bar x'$ of (\ref{eq:subspace3}) is finite,
and this implies the proposition in this case.

Now suppose that $c_1=0$ in (\ref{eq:lin}). One of $c_2,\ldots,c_s$ is non-zero, and  for simplicity,
we assume that $c_s\ne 0$. Then combining (\ref{eq:subspace}) with (\ref{eq:lin}),
we obtain that 
\begin{equation}\label{eq:subspace4}
\left|b_1 +\sum_{j=2}^{s-1} (b_j-c_jb_sc_s^{-1}) x_j\right|_v<\hbox{H}(\bar x)^{-\epsilon}.
\end{equation}
Given a solution $\bar x$ of (\ref{eq:subspace4}), we pick a minimal 
$J\subset \{2.\ldots, s-1\}$ such that
\begin{equation}\label{eq:subspace5}
\left|b_1 +\sum_{j\in J} (b_j-c_jb_sc_s^{-1}) x_j\right|_v<\hbox{H}(\bar x)^{-\epsilon},
\end{equation}
and no proper subsum of $b_1 +\sum_{j\in J} (b_j-c_jb_sc_s^{-1}) x_j$ vanishes.
Let $\bar x'=(x_j: j\in J)$. Since $\hbox{H}(\bar x')\le \hbox{H}(\bar x)$,
it follows from the inductive hypothesis that $\bar x'$ belongs to a fixed finite set.
This proves that for every solution $\bar x$ of (\ref{eq:subspace})
there exists $j_0=2,\ldots, s$ such that $x_{j_0}$ belongs to a fixed finite set.
Now we can finish the argument as in the previous paragraph, and this completes
the proof of the proposition.
\QED

\subsection{Higher order mixing under irreducibility condition}\label{sec:mult_irr}
 We define the subspace $W $ in $\mathcal{L}(G)$, the set of characters $\mathcal{X}_0$ and  
the eigenvectors $w_\chi$ with $\chi\in \mathcal{X}_0$
as in Section \ref{sec:3_irred}. 

In this section we assume that $D\pi(W)$ is not contained in any proper rational subspace.
Let $\{w_1,\ldots,w_k\}$ be a fixed basis of $W$. 
Consider a box map
$$
\iota:B\to \mathcal{L}(G): t\mapsto \sum_{i=1}^k t_i w_i
$$
where $B=[0,T_1]\times \cdots \times [0,T_k]$.

\begin{lemma}\label{l:m_box}
Let $f_1,\ldots,f_s\in C(X)$, $u_1,\ldots,u_s\in\mathcal{L}(G)$,
$\beta_1,\ldots,\beta_s$ be automorphisms of $G$ such that $\beta_i=id$ of $G/G'$,
$z_1,\ldots z_s\in \Z^l$, $v_1,\ldots, v_s\in\mathcal{L}(G)$, and $x_1,\ldots,x_s\in X$. Then 
$$
\frac{1}{|B|} \int_B \left( \prod_{i=1}^s f_i(\exp(u_i)\beta_i(\alpha(z_i)(\exp(v_i+\iota(t))))x_i) \right)\, dt
= \prod_{i=1}^s \left(\int_X f_i\,d\mu\right)+o(1)
$$
as $\min\{\|z_i\|, \|z_i-z_j\|: i\ne j\}\to \infty$.
Moreover, the convergence is uniform over $u_1,\ldots, u_s$,
$\beta_1,\ldots,\beta_s$, $v_1,\ldots,v_s$, $x_1,\ldots x_s$,
and functions $f_1,\ldots, f_s$ with $\|f_i\|_{C^\theta}\ll 1$.
\end{lemma}

\proof
It is sufficient to prove the claim when $f_1,\ldots,f_s$ belong to a dense family of functions in $C(X)$.
Hence, without loss of generality, we may assume that the functions are H\"older with exponent $\theta$.

To prove the lemma, we apply Theorem \ref{th:green-tao} to the product nilmanifold $X^s=G^s/\Lambda^s$.
Suppose that the claim of the lemma fails. Then there exist $\delta\in (0,\delta_0)$ 
and sequences $f^{(n)}_1,\ldots,f^{(n)}_s\in C^\theta(X)$, $u^{(n)}_1,\ldots,u^{(n)}_s\in\mathcal{L}(G)$,
$\beta^{(n)}_1,\ldots,\beta^{(n)}_s$ satisfying $\beta^{(n)}_i=id$ on $G/G'$,
$z_1^{(n)},\ldots, z_s^{(n)}\in \Z^l$ satisfying
$$
\min\left\{\|z^{(n)}_i\|, \|z^{(n)}_i-z^{(n)}_j\|: i\ne j\right\}\to \infty\quad\hbox{as $n\to\infty$,}
$$
$v_1^{(n)},\ldots, v_s^{(n)}\in\mathcal{L}(G)$, 
and $x^{(n)}_1,\ldots,x^{(n)}_s\in X$ such that 
\begin{align*}
&\left|\frac{1}{|B|} \int_{B} \left( \prod_{i=1}^s
    f_i^{(n)}(\exp(u^{(n)}_i)\beta^{(n)}_i(\alpha(z^{(n)}_i)(\exp(v_i^{(n)}+\iota(t))))x_i^{(n)}) \right)\, dt -
  \prod_{i=1}^s \left(\int_X f_i^{(n)}\,d\mu\right)\right|\\
 >&\delta \prod_{i=1}^s\|f_i^{(n)}\|_{C^\theta}.
\end{align*}
We set $f^{(n)}=f^{(n)}_1\otimes \cdots \otimes f^{(n)}_s:X^s\to \R$, $u^{(n)}=(u^{(n)}_1,\ldots,u^{(n)}_s)\in\mathcal{L}(G)^s$,
$$
\iota^{(n)}:B\to\mathcal{L}(G)^s: t\mapsto \left(D\beta^{(n)}_1 D\alpha(z^{(n)}_1)(v^{(n)}_1+\iota(t)),\ldots,
D\beta^{(n)}_s D\alpha(z^{(n)}_s)(v^{(n)}_s+\iota(t))\right),
$$
and $x^{(n)}=(x^{(n)}_1,\ldots,x^{(n)}_s)\in X^s$.
Then
$$
\|f^{(n)}\|_{C^\theta}\ll \prod_{i=1}^s\|f_i^{(n)}\|_{C^\theta},
$$
and
\begin{align*}
\left|\frac{1}{|B|} \int_{B} 
    f^{(n)}(\exp(u^{(n)})\exp(\iota^{(n)}(t))x^{(n)})\, dt -
  \int_X f^{(n)}\,d\mu^{\otimes s}\right|  \gg\delta \|f^{(n)}\|_{C^\theta}.
\end{align*}
It follows from Theorem \ref{th:green-tao} that there exists $(a_1^{(n)},\ldots,a_s^{(n)})\in (\Z^d)^s\backslash
\{0\}$ such that 
\begin{equation}\label{eq:bound_1}
\|a_1^{(n)}\|,\ldots, \|a_s^{(n)}\|\ll \delta^{-L_1}\ll 1,
\end{equation}
and 
\begin{equation}\label{eq:bound_2}
\left|\sum_{j=1}^s \left<a^{(n)}_j, D\pi D\beta^{(n)}_j D\alpha(z_j^{(n)})(w_i) \right> \right|\ll
\delta^{-L_2}/T_i\ll 1
\quad\hbox{for all $i=1,\ldots,k$.}
\end{equation}

Since $\beta_j^{(n)}=id$ on $G/G'$, we have $D\pi D\beta^{(n)}_j=D\pi$.
We rewrite (\ref{eq:bound_2}) in terms of vectors $w_\chi$, $\chi\in \mathcal{X}_0$ that
satisfy $D\alpha(z)w_\chi=\chi(z)w_\chi$ for $z\in \Z^l$. Since each $w_\chi$ can be written as
a linear combination of $w_i$'s, it follows from (\ref{eq:bound_2}) that 
\begin{equation}\label{eq:bound_21}
\left|\sum_{j=1}^s \chi (z_j^{(n)})\left<a^{(n)}_j, D\pi (w_\chi) \right> \right| \ll 1
\quad\hbox{for all $\chi\in\mathcal{X}_0$.}
\end{equation}
We observe that because of (\ref{eq:bound_1}), the tuple $(a_1^{(n)},\ldots,a_s^{(n)})$ varies
over a finite set. Hence, passing to subsequence, we may assume that (\ref{eq:bound_21}) 
holds for a fixed tuple $(a_1,\ldots,a_s)\in (\Z^d)^s\backslash \{0\}$.
After changing indices, we may assume that $a_j\ne 0$ for $j=1,\ldots s'$ and $a_j=0$ for $j>s'$.
We note that this implies that 
$$
b_j:=\left<a_j,D\pi(w_\chi)\right>\ne 0\quad\hbox{for all $j=1,\ldots s'$ and $\chi\in\mathcal{X}_0$.}
$$ 
Indeed, if $\left<a_j,D\pi(w_\chi)\right>= 0$ for some $j$ and $\chi$, then
taking Galois conjugates we obtain that  $\left<a_j,D\pi(w_\chi)\right>= 0$
for all $\chi\in\mathcal{X}_0$. This implies that $D\pi(W)$ is contained in a proper rational subspace
and contradicts the irreducibility assumption. 

We may cancel vanishing subsums from (\ref{eq:bound_21}), and 
passing to a subsequence, we may assume that no proper subsum in (\ref{eq:bound_21}) vanishes.

Passing to a subsequence and changing indices, we may also assume that
\begin{equation}\label{eq:max}
\max\left\{\|z^{(n)}_j\|: 1\le j\le s'\right\}=\|z_1^{(n)}\|.
\end{equation}
By Lemma \ref{l:expansion}, there exists fixed $c>0$ such that
\begin{equation}\label{eq:expa}
\max_{\chi\in\mathcal{X}_0} |\chi(z)|\ge \exp(c\|z\|),\quad z\in \Z^l.
\end{equation}
Hence, passing to a subsequence, we may assume that 
$$
|\chi_0(z_1^{(n)})|\ge \exp(c\|z_1^{(n)}\|)
$$
holds with a fixed $\chi_0\in\mathcal{X}_0$. For this $\chi_0$, (\ref{eq:bound_21}) gives
\begin{equation}\label{eq:bound_22}
\left|b_0+\sum_{j=1}^{s'} b_j x^{(n)}_j \right| \ll \exp(-c\|z_1^{(n)}\|),
\end{equation}
where $x^{(n)}_j:=\chi_0 (z_j^{(n)}-z_1^{(n)})$.
It is clear that $b_j$'s and $x_j^{(n)}$ are $S$-units in a fixed number field,
and to derive a contradiction, we apply the estimate of Proposition \ref{p:subspace}.
We observe that there exists $c_v>1$, $v\in S$, such that
$$
|\chi_0(z)|_v\le \exp(c_v\|z\|),\quad z\in \Z^l.
$$
Hence,
$$
|x_j^{(n)}|_v=|\chi_0 (z_j^{(n)}-z_1^{(n)})|_v\le \exp(c_v \|z_j^{(n)}-z_1^{(n)}\|),
$$
and by (\ref{eq:max}),
\begin{align*}
\hbox{H}(\bar x^{(n)})&=\prod_{v\in S}  \max (1, |x_1^{(n)}|_v,\ldots,|x_{s'}^{(n)}|_v)\\
&\le \exp\left(\left(\sum_{v\in S} c_v\right)\max_{1\le j\le s'} \|z_j^{(n)}-z_1^{(n)}\| \right)\\
&\le \exp\left(2 \left(\sum_{v\in S} c_v\right)\|z_1^{(n)}\| \right).
\end{align*}
It follows from (\ref{eq:bound_22}) that 
\begin{equation}\label{eq:bound_23}
\left|b_0+\sum_{j=1}^{s'} b_j x^{(n)}_j \right| \ll \hbox{H}(\bar x^{(n)})^{-\epsilon}.
\end{equation}
with fixed $\epsilon>0$. According to our construction, no proper subsum in (\ref{eq:bound_23}) vanishes.
Hence, it follows from Proposition \ref{p:subspace} that $x^{(n)}=\chi_0 (z_j^{(n)}-z_1^{(n)})$ runs over a finite
set. Since all elements in $\mathcal{X}_0$ are conjugate under the Galois action,
it follows that $\chi (z_j^{(n)}-z_1^{(n)})$ with $\chi\in\mathcal{X}_0$ also runs over finite set.
In particular, 
$$
\max_{\chi\in\mathcal{X}_0} |\chi (z_j^{(n)}-z_1^{(n)})|\ll 1.
$$
On the other hand, by (\ref{eq:expa}),
$$
\max_{\chi\in\mathcal{X}_0} |\chi (z_j^{(n)}-z_1^{(n)})|\to \infty.
$$
This contradiction proves the lemma.   \QED

Now we prove Theorem \ref{th:multmixing} under the irreducibility condition.
Without loss of generality, we may assume that $z_0=0$.
We fix a fundamental domain $F\subset G$ for $G/\Lambda$ and set $E=\exp^{-1}(F)$.
We may arrange that $E$ is bounded and has piecewise smooth boundary.
Then
$$
\int_X f_0(x)\left(\prod_{i=1}^s f_i(\alpha(z_i)(x))\right)\, d\mu(x)
=\int_E f_0(\exp(u)\Lambda)\left(\prod_{i=1}^s f_i(\exp(D\alpha(z_i)u)\Lambda)\right)\, du,
$$
where $du$ denotes a suitably normalised Lebesgue measure on $\mathcal{L}(G)$.
We choose a basis of $\mathcal{L}(G)$ that contains the fixed basis $\{w_i\}$ of $W$
and tessellate $\mathcal{L}(G)$ by cubes $C$ of size $\epsilon$ with respect to this basis.
Then
\begin{equation}\label{eq:boundary}
\left|E-\bigcup_{C\subset E} C\right|\ll \epsilon.
\end{equation}
For all $u_1,u_2\in C$,
\begin{equation}\label{eq:boundary2}
|f_0(\exp(u_1)\Lambda)-f_0(\exp(u_2)\Lambda)|\ll \epsilon^\theta .
\end{equation}
Here and later in the argument the implied constants may depend on the H\"older norms of $f_0,\ldots,f_s$.
For every cube $C$, we pick a point $u_C\in C$. Then it follows from (\ref{eq:boundary}) and
(\ref{eq:boundary2}) that
\begin{align}\label{eq:expand}
&\int_E f_0(\exp(u)\Lambda)\left(\prod_{i=1}^s f_i(\exp(D\alpha(z_i)u)\Lambda)\right)\, du\\
=& \sum_{C\subset E} \int_C f_0(\exp(u)\Lambda)\left(\prod_{i=1}^s
  f_i(\exp(D\alpha(z_i)u)\Lambda)\right)\, du +O(\epsilon) \nonumber\\
=& \sum_{C\subset E} f_0(\exp(u_C)\Lambda) \int_C \left(\prod_{i=1}^s
  f_i(\exp(D\alpha(z_i)u)\Lambda)\right)\, du +O(\epsilon^\theta).\nonumber
\end{align}
Each cube $C$ in the above sum can be written as $C=B'+B$ where $B$ is a cube in $W$ and
$B'$ is a cube in the complementary subspace. 

Let
$$
N=N(z_1,\ldots, z_s):=\min(\|z_i\|,\|z_i-z_j\|:\, i\ne j).
$$
It follows from Lemma \ref{l:m_box} that
\begin{equation}\label{eq:lemm}
\frac{1}{|B|} \int_B \left(\prod_{i=1}^s f_i(\exp(v_i+D\alpha(z_i)(b))\Lambda)  \right)\, db
\to \prod_{i=1}^s \left(\int_X f_i\, d\mu\right)
\end{equation}
as $N\to\infty$, uniformly over $v_1,\ldots v_s\in \mathcal{L}(G)$ and the cubes $B$
(note that all the cubes are translates of a fixed cube).
Hence, it follows that
\begin{align*}
&\frac{1}{|C|} \int_C \left(\prod_{i=1}^s f_i(\exp(D\alpha(z_i)(u))\Lambda) \right)\, du\\
=&\frac{1}{|B'||B|} \int_{B'}\int_B \left(\prod_{i=1}^s
  f_i(\exp(D\alpha(z_i)(b')+D\alpha(z_i)(b))\Lambda) \right)\, dbdb'\\
\to& \prod_{i=1}^s \left(\int_X f_i\,d\mu\right)
\end{align*}
as $N\to \infty$. Combining this with (\ref{eq:expand}), we deduce that
\begin{align*}
&\int_E f_0(\exp(u)\Lambda)\left(\prod_{i=1}^s f_i(\exp(D\alpha(z_i)u)\Lambda)\right)\, du\\
=& 
\left(\sum_{C\subset E} f_0(\exp(u_C)\Lambda) |C| \right)\prod_{i=1}^s \left(\int_X f_i\,d\mu\right)
+ \left(\sum_{C\subset E} |C|\right) o(1)+O(\epsilon^\theta),
\end{align*}
where clearly $\sum_{C\subset E} |C|=O(1)$.
Moreover, using (\ref{eq:boundary2}) and (\ref{eq:boundary}), we deduce that
\begin{align*}
\sum_{C\subset E} f_0(\exp(u_C)\Lambda) |C|&= \sum_{C\subset E} \left(\int_C f_0(\exp(u)\Lambda)\, du +O(|C|\epsilon^\theta)\right)\\
&=\int_E f_0(\exp(u)\Lambda)\, du+O(\epsilon^\theta)\\
&=\int_X f_0\, d\mu+O(\epsilon^\theta).
\end{align*}
This implies that
\begin{align*}
\int_X f_0(x)\left(\prod_{i=1}^s f_i(\alpha(z_i)(x))\right)\, d\mu(x)
&=\int_E f_0(\exp(u)\Lambda)\left(\prod_{i=1}^s f_i(\exp(D\alpha(z_i)u)\Lambda)\right)\, du\\
&=\prod_{i=0}^s \left(\int_X f_i\,d\mu\right)+o(1)+O(\epsilon^\theta).
\end{align*}
as $N\to\infty$. This proves Theorem
\ref{th:multmixing}
under the irreducibility condition.
It is clear from the proof that convergence is uniform provided that
$\|f_0\|_{C^\theta},\ldots, \|f_s\|_{C^\theta}\ll 1$.

\subsection{Higher order mixing in general}\label{sec:mult_gen}
We will apply an inductive argument which uses the result of Section \ref{sec:mult_irr}
as a base case. In fact, we note that the argument in Section \ref{sec:mult_irr}
implies that 
\begin{equation}\label{eq:conv_gen}
\int_X f_0(x) \left(\prod_{i=1}^s f_i(h_i\beta_i(\alpha(z_i)(x)))\right)\, d\mu(x)
= \prod_{i=0}^s \left(\int_X f_i\, d\mu\right)+o(1)
\end{equation}
as $N=N(z_1,\ldots,z_s)\to \infty$, uniformly over 
functions $f_0,\ldots,f_s$ with $\|f_i\|_{C^\theta}\ll 1$,
$h_1,\ldots, h_s\in G$ and
automorphisms $\beta_1,\ldots,\beta_s$ of $G$ that preserve $\Lambda$ and act trivially on $G/G'$.
Indeed, Lemma \ref{l:m_box} implies that in (\ref{eq:lemm}) we more generally have
\begin{align*}
\frac{1}{|B|} \int_B \left(\prod_{i=1}^s f_i(h_i\beta_i(\exp(v_i+D\alpha(z_i)(b)))\Lambda)  \right)\, db
= \prod_{i=1}^s \left(\int_X f_i\, d\mu\right)+o(1)
\end{align*}
as $N\to\infty$, uniformly over $f_0,\ldots,f_s$ with $\|f_i\|_{C^\theta}\ll 1$,
$h_1,\ldots, h_s\in G$, automorphisms $\beta_1,\ldots,\beta_s$,
and $v_1,\ldots v_s\in \mathcal{L}(G)$. Then the rest of the argument carries over
 and implies (\ref{eq:conv_gen}).

Let $W$ be the Lie subalgebra of $\mathcal{L}(G)$ introduced in Section \ref{sec:3_irred}.
By \cite[Ch.~5, Sec.~5]{Starkov}, there exists a closed connected normal subgroup $M$ of $G$ such that
$M/(M\cap \Lambda)$ is compact, and 
$$
\overline{\exp(W)g\Lambda}=Mg\Lambda\quad\hbox{for almost every $g\in G$.}
$$
Since we may replace the lattice $\Lambda$ by its conjugate, we assume that 
$$
\overline{\exp(W)\Lambda}=M\Lambda.
$$
We note that the group $M$ satisfies the following properties:
\begin{enumerate}
\item[(i)] $M$ is $\alpha(\Z^l)$-invariant,
\item[(ii)] $D\pi(W)$ is not contained in a proper rational subspace of $\mathcal{L}(M/M')$,
\item[(iii)] $[G,M]\subset M'$.
\end{enumerate}
Properties (i)--(iii) can be verified exactly as in the proof of \cite[Lem.~3.4]{GS}.

To apply induction, we observe that the nilmanifold $X=G/\Lambda$ fibers over the nilmanifold of
$Y=G/(M\Lambda)$ with fibers isomorphic to $R=M\Lambda/\Lambda\simeq M/(M\cap\Lambda)$.
The invariant measure $\mu$ on $X$ decomposes as
$$
\int_X f\, d\mu=\int_Y \int_R f(yr)\, d\mu_R(r)d\mu_Y(y),\quad f\in C(X),
$$
where $\mu_Y$ and $\mu_R$ denote the normalised invariant measures on $Y$ and $R$ respectively. 
It follows from (i) that the fibration $X\to Y$ is $\alpha(\Z^l)$-equivariant.
Hence, we obtain
\begin{align}\label{eq:dd1}
&\int_X f_0(x)\left(\prod_{i=1}^s f_i(\alpha(z_i)(x))\right)\, d\mu(x)\\
=& \int_Y\left(\int_R  f_0(yr)\left(\prod_{i=1}^s f_i(\alpha(z_i)(y)\alpha(z_i)(r))\right)d\mu_R(r)
\right) d\mu_Y(y)\nonumber\\
=& \int_F\left(\int_R  f_0(gr)\left(\prod_{i=1}^s f_i(\alpha(z_i)(g)\alpha(z_i)(r))\right)d\mu_R(r)
\right) dm_F(g),\nonumber
\end{align}
where $F\subset G$ is a bounded fundamental set for $G/(M\Lambda)$, and $m_F$ is the measure on $F$
induced by $\mu_Y$. 
We write
$$
\alpha(z_i)(g)=a_im_i\lambda_i\quad \hbox{with $a_i\in F$, $m_i\in M$, $\lambda_i\in \Lambda$},\quad\quad
i=1,\ldots,s.
$$
Then 
\begin{align*}
&\int_R  f_0(gr)\left(\prod_{i=1}^s f_i(\alpha(z_i)(g)\alpha(z_i)(r))\right)d\mu_R(r)\\
=&\int_R  f_0(gr)\left(\prod_{i=1}^s f_i(a_im_i\beta_i(\alpha(z_i)(r)))\right)d\mu_R(r),
\end{align*}
where $\beta_i$'s are the transformations of $S$ induced by the automorphisms $m\mapsto
\lambda_im\lambda_i^{-1}$ of $M$. By (ii), the subspace $W\subset \mathcal{L}(M)$ satisfies
the irreducibility assumption of Section \ref{sec:mult_irr}, and by (iii), the automorphisms
$\beta_i$ act trivially on $M/M'$. Let
$$
\phi_0(r):=f_0(gr)\quad\hbox{and}\quad \phi_i(r):=f_i(a_ir),\quad i=1,\ldots,s.
$$
Since $F\subset G$ is bounded, we have $\|\phi_i\|_{C^\theta}\ll \|f_i\|_{C^\theta}\ll 1$.
Hence, it follows from (\ref{eq:conv_gen}) that 
\begin{align*}
\int_R  \phi_0(r)\left(\prod_{i=1}^s \phi_i(m_i\beta_i(\alpha(z_i)(r)))\right)d\mu_R(r)
= 
\prod_{i=0}^s \left(\int_R \phi_i\, d\mu_R\right)+o(1)
\end{align*}
as $N\to \infty$, uniformly over $g\in F$, $m_1,\ldots,m_s\in M$,
and the automorphisms $\beta_1,\ldots, \beta_s$.
Since $a_i M\Lambda=\alpha(z_i)(g)M\Lambda$,  this implies that
\begin{align}\label{eq:dd2}
\int_R  f_0(gr)\left(\prod_{i=1}^s f_i(\alpha(z_i)(g)\alpha(z_i)(r))\right)d\mu_R(r)
= \prod_{i=0}^s \left(\int_R f_i(\alpha(z_i)(g)r)\, d\mu_R(r)\right)+o(1)
\end{align}
as $N\to\infty$, uniformly over $g\in F$. Let $\bar f_i$ be the function on $Y$ defined by
$$
y\mapsto \int_R f_i(yr)\, d\mu_R(r).
$$
Combining (\ref{eq:dd1}) and (\ref{eq:dd2}), we deduce that
$$
\int_X f_0(x) \left(\prod_{i=1}^s f_i(\alpha(z_i)(x))\right)\, d\mu(x)
=\int_Y \bar f_0(y) \left(\prod_{i=1}^s \bar f_i(\alpha(z_i)(y))\right)\, d\mu_Y(y)+o(1)
$$
as $N\to\infty$. Finally, it follows by induction by $\dim(X)$ that
\begin{align*}
\int_Y \bar f_0(y)\left(\prod_{i=1}^s \bar f_i(\alpha(z_i)(y))\right)\, d\mu_Y(y)
&=\prod_{i=0}^s \left( \int_Y \bar f_i\, d\mu_Y\right)+o(1)\\
&=\prod_{i=0}^s \left( \int_X f_i\, d\mu\right)+o(1).
\end{align*}
The above argument implies uniform convergence provided that
$\|f_0\|_{C^\theta},\ldots, \|f_s\|_{C^\theta}\ll 1$.
This completes the proof of Theorem \ref{th:multmixing}.

\section{Exponential mixing of shapes}   \label{Exponential multiple mixing}

While we have proved exponential 2-mixing and 3-mixing for $\Z^l$-actions by automorphisms on
nilmanifolds, we do not know if exponential mixing of higher orders
holds for them in general.  This would require to establish
a quantitative version of Proposition \ref{p:subspace}
which seems to be out of reach of available number-theoretic methods.  
Nonetheless, we prove a weak form of exponential mixing where the error term is controlled
by
\[
N_{*} (z_0,\ldots,z_s):=\min_{\chi\in \mathcal{X}(\alpha)} \:\min_{|\chi(z_i-z_j)|\ge 1} \{|\chi(z_i-z_j)| : i \ne j\} 
\]
with notation as in Section \ref{sec:notation}.
This, in particular, implies exponential mixing for Anosov shapes --- Theorem \ref{th:shape}.



\begin{Theorem} \label{th:multmixing_2}
Let $\alpha:\Z^l\to \hbox{Aut}(X)$ be an action on a compact nilmanifold $X$
such that every $\alpha(z)$, $z\ne 0$, is ergodic.
Then there exists $\eta=\eta(\theta)>0$ such that for every $f_0,\ldots,f_s\in C^\theta(X)$ and $z_0,\ldots,z_s\in \mathbb{Z}^l$,
\begin{equation}  \label{eq:mult-irr}
\int_X \left( \prod_{i=0}^s f_i(\alpha (z_i)x)\right)\, d\mu(x)
= \prod_{i=0}^s \left(\int_X f_i\, d\mu\right)+
O\left( N_{*}(z_0,\ldots,z_s)^{-\eta} \prod_{i=0}^s \|f_i\|_{C^\theta} \right).
\end{equation}
\end{Theorem}
\vspace{.5em}

\proof  
We adapt the method of the proof of \cite[Th.~1.2]{GS} from our previous paper.
Since the proof is quite similar to the argument in this paper,
we will only give an outline. 

We take a character $\chi\in\mathcal{X}(\alpha)$
and the corresponding eigenvector $w\in\mathcal{L}(G)$. If $\chi$ is real, we denote by
$W$ the corresponding one-dimensional eigenspace. Otherwise, we denote by $W$ the two-dimensional
subspace $\left<w,\bar w\right>\cap \mathcal{L}(G)$. Then $D\alpha(z)|_W=r(z)\cdot \omega(z)$
where $r(z)=|\chi(z)|$ and $\omega(z)$ is a rotation. We assume in addition that $\chi^2\notin
\mathcal{X}(\alpha)$. Then $W$ is closed under the Lie bracket.

We first treat the irreducible case: namely, when $D\pi(W)$ is not contained in a proper rational subspace.
Without loss of generality, $N_*>1$. Then for all $i\ne j$,
$|\chi(z_i-z_j)| \ne  1 $ and after changing indexes, $| \chi (z_0)| < | \chi (z_1)| < \ldots < | \chi (z_s) |$. 
We may also assume that $z_0=0$. We fix a basis of $\mathcal{L}(G)$ which contains a basis of $W$
and tessellate $\mathcal{L}(G)$ by cubes of size $\epsilon$ with respect to this basis.
Then the integral
$$
\int_X f_0(x)\left( \prod_{i=1}^s f_i(\alpha (z_i)x)\right)\, d\mu(x)
$$
can be approximated by a sum of the integrals
\begin{equation*}
\int_C \left( \prod_{i=1}^s f_i(\exp(D\alpha (z_i)u)\Lambda)\right)\,du
\end{equation*}
with the error of size $O(\epsilon^\theta \prod_{i=1}^s \|f\|_{C^\theta})$.
Since $C=B'+B$ where $B$ is a cube in $W$ and $B'$ is a cube in the complementary subspace,
the above integral can be written as
$$
\int_{B'}\int_B \left( \prod_{i=1}^s f_i(\exp(D\alpha (z_i)b'+D\alpha (z_i)b)\Lambda)\right)\,dbdb'.
$$
For every cube $B$, we take a box map $\iota_B:[0,\epsilon]^{\dim(W)}\to B$ that parametrises $B$. Then
by \cite[Prop.~4.2]{GS}, there exists $\kappa>0$ such that
\begin{align*}
&\frac{1}{|B|} \int_B \left( \prod_{i=1}^s f_i(\exp(v_i+D\alpha (z_i)b)\Lambda)\right)\,db\\
=& \epsilon^{-\dim(W)} \int_{[0,\epsilon]^{\dim(W)}} \left( \prod_{i=1}^s
  f_i(\exp(v_i+r(z_i)\omega(z_i)\iota_B(t))\Lambda)\right)\,dt\\
=& \prod_{i=1}^s \left(\int_X f_i\, d\mu\right)+O\left(\sigma^{-\kappa}\prod_{i=1}^s \|f\|_{C^\theta} \right)
\end{align*}
uniformly over $v_1,\ldots,v_s\in\mathcal{L}(G)$, where
$\sigma=\min\{\epsilon r(z_1), r(z_s)/r(z_{s-1}),\ldots, r(z_2)/r(z_1)\}$.
We note that the Diophantine condition for the box map $t\mapsto \omega(z_i)\iota_B(t)$,
which is required in \cite[Prop.~4.2]{GS}, is satisfied because $W$ is spanned by vectors with
algebraic coordinates and \cite[Th.~7.3.2]{BG} applies.
Since this estimate is uniform over $v_i$'s, it follows that
$$
\frac{1}{|C|}\int_C \left( \prod_{i=1}^s f_i(\exp(D\alpha (z_i)u)\Lambda)\right)\,du
= \prod_{i=1}^s \left(\int_X f_i\, d\mu\right)+O\left(\sigma^{-\kappa}\prod_{i=1}^s \|f\|_{C^\theta} \right),
$$
and we deduce that
$$
\int_X f_0(x)\left( \prod_{i=1}^s f_i(\alpha (z_i)x)\right)\, d\mu(x)
= \prod_{i=1}^s \left(\int_X f_i\, d\mu\right)+O\left((\epsilon^\theta+\sigma^{-\kappa})\prod_{i=1}^s \|f\|_{C^\theta} \right).
$$
We refer to the proof of \cite[Th.~1.2]{GS} for details.
Choosing $\epsilon=r(z_1)^{-\frac{1}{2}}$ implies the claim of the theorem in the irreducible case.

To give a proof in general, we use induction on $\dim(X)$.
This argument is very similar to  Section \ref{sec:3_gen}. If $\overline{\exp(W)\Lambda}\ne X$, 
we consider the $\alpha$-equivariant fibration $X\to Y$ defined by the closure.
The above argument implies that
$$
\int_X f_0(x)\left( \prod_{i=1}^s f_i(h_i\beta_i(\alpha (z_i)x))\right)\, d\mu(x)
= \prod_{i=1}^s \left(\int_X f_i\, d\mu\right)+O\left(N_*^{-\eta}\prod_{i=1}^s \|f\|_{C^\theta} \right)
$$
uniformly over $h_1,\ldots,h_s\in G$ and automorphisms $\beta_1,\ldots,\beta_s\in\hbox{Aut}(X)$
that act trivially on $G/G'$.  (In fact, this uniformity was part of  \cite[Prop.~4.2]{GS}.)
Then as in Section \ref{sec:3_gen},
\begin{align*}
\int_X f_0(x)\left( \prod_{i=1}^s f_i(\alpha (z_i)x)\right)\, d\mu(x)
=\int_Y \bar f_0(x)\left( \prod_{i=1}^s \bar f_i(\alpha (z_i)x)\right)\, d\mu_Y(x)
+O\left(N_*^{-\eta}\prod_{i=1}^s \|f\|_{C^\theta} \right)
\end{align*}
where $\bar f_0,\ldots,\bar f_s\in C^\theta(Y)$. Now the claim follows by induction on dimension.
\QED

\begin{remark}
{\rm 
In the irreducible case of the above proof, we can replace $N_{*} (z_0,\ldots,z_s)$
\[
N'_{*} (z_0,\ldots,z_s):=\max_{\chi\in \mathcal{X}'(\alpha)} \:\min_{|\chi(z_i-z_j)|\ge 1} \{|\chi(z_i-z_j)| : i \ne j\},
\]
which provides a better estimate.
}
\end{remark}

{\em Proof of Theorem \ref{th:shape}}\hspace{.1em}: We note that if in Theorem \ref{th:multmixing_2} we assume that $\alpha(z_i-z_j)$
are Anosov for all $i\ne j$, then $N_*(z_0,\ldots,z_s)>1$ and 
$N_*(nz_0,\ldots,nz_s)=N_*(z_0,\ldots,z_s)^n$. Hence, Theorem \ref{th:shape} follows directly 
from Theorem \ref{th:multmixing_2}.

\section{Cocycle rigidity}

We now apply exponential 2-mixing to prove smooth cocycle  rigidity for genuinely higher-rank abelian
actions by automorphisms of nilmanifolds --- Theorem \ref{th:cocyle}.
The proof is based on the ``higher-rank trick'' from \cite{KS}.

Let $\Z^l\to \hbox{Aut}(X)$ be an action on a compact nilmanifold $X$
and $c:\Z^l\times X\to \R$ a cocycle.
Assume that there is a rank-two subgroup $\left<a,b\right>$ of $\Z^l$ such that 
its every nonzero element acts ergodically on $X$. 

First, we note that the map $c_0: \Z^l \mapsto \R$ defined by 
$z \mapsto \int_X c (z,x)\, d\mu(x)$
is a homomorphism by the cocycle property.  Then $c -c_0$ is also a cocycle,
 and it will suffice to prove the theorem for $c -c_0$.  Thus, we will assume that all functions $c
 (z,\cdot)$ for $z \in \Z^l$ satisfy $\int _X c (z,x)\, d\mu(x) =0$.

Let $f(x)=c(a,x)$.
We shall show that there exists $\phi\in L^2(X)$ such that
\begin{equation}\label{eq:cob}
f=\phi\circ a-\phi.
\end{equation}
We will apply our previous results  \cite[Sec.~6]{GS}.  By \cite[Theorem 6.1]{GS}, it suffices to show that 
$$
 \sigma ^2 := \int _X f^2\,  d\mu  +  2  \sum _{i=1}  ^{\infty}   \int _X (  f \circ a ^i )  f\,  d\mu 
= \sum _{i=-\infty}  ^{\infty}   \left< f \circ a ^i ,  f\right>=0.
$$
We note that the assumption of this theorem is verified in \cite[Sec.~6]{GS} using exponential mixing of $a$.
Let $h(x)=c(b^j,x)$. It follows from the cocycle property that
$$
f\circ b^j-f=h\circ a-h.
$$
Hence,
\begin{align*}
\sum_{i=-n}^n (f\circ a^ib^j-f\circ a^i)
=\sum_{i=-n}^n (h\circ a^{i+1}-h\circ a^i)=h\circ a^{n+1}-h\circ a^{-n},
\end{align*}
and it follows from exponential mixing that 
$$
\sigma ^2= \sum _{i=-\infty}  ^{\infty}   \left< f \circ a ^ib^j ,  f\right>\quad
\hbox{for every $j\in\Z$.}
$$
 On the other hand, by the exponential mixing for the group $\left<a,b\right>$
established in Theorem \ref{th:3_mixing},
$$
\sum _{i,j\in\Z}   \left< f \circ a ^ib^j ,  f\right><\infty.
$$
This implies that $\sigma^2=0$ and proves (\ref{eq:cob}).

Now using the cocycle regularity result \cite[Theorem 7.1]{GS} established in our previous paper,
we deduce that (\ref{eq:cob}) also has a $C^\infty$ solution.
Finally, since $a$ acts ergodically, it follows from \cite[Lem.~4.1]{KS} that
$$
c(z,\cdot)=\phi\circ z-\phi\quad\hbox{for all $z\in \Z^l$.}
$$
This completes the proof of Theorem \ref{th:cocyle}.




\bibliographystyle{alpha}

\end{document}